\documentclass{amsart}
\usepackage{a4wide}
\usepackage[T1]{fontenc}
\usepackage{amsmath,amssymb,amsthm,stmaryrd}
\usepackage{mathrsfs}
\usepackage{latexsym}
\usepackage{esint}
\usepackage{color}

\newcommand{\ud}[0]{\,\mathrm{d}}

\newcommand{\eps}[0]{\varepsilon}

\newcommand{\abs}[1]{|#1|}
\newcommand{\Babs}[1]{\Big|#1\Big|}

\newcommand{\BNorm}[2]{\Big\|#1\Big\|_{#2}}

\renewcommand{\Im}[0]{\operatorname{Im}}

\newcommand\R{\mathbb{R}}
\newcommand\C{\mathbb{C}}

\newcommand\Z{\mathbb{Z}}

\newcommand\Q{\mathbb{Q}}

\newcommand{\prob}[0]{\mathbb{P}}

\swapnumbers

\numberwithin{equation}{section}

\makeatletter
  \let\c@subsection\c@equation
\makeatother

\theoremstyle{plain}
\newtheorem{theorem}[equation]{Theorem}
\newtheorem{proposition}[equation]{Proposition}
\newtheorem{corollary}[equation]{Corollary}
\newtheorem{lemma}[equation]{Lemma}

\theoremstyle{definition}

\theoremstyle{remark}

\usepackage[english]{babel}

\usepackage[letterpaper,top=2cm,bottom=2cm,left=3cm,right=3cm,marginparwidth=1.75cm]{geometry}

\usepackage{amsmath}
\usepackage{graphicx}
\usepackage[colorlinks=true, allcolors=blue]{hyperref}

\title{\bf Guido Weiss: a few memories of a friend and an influential mathematician
}

\begin{document}
\author{Pascal Auscher}
\address{Universit\'e Paris-Saclay, CNRS, Laboratoire de Math\'{e}matiques d'Orsay, 91405 Orsay, France}
\email{pascal.auscher@universite-paris-saclay.fr} 
\author{Aline Bonami}
\address{ Institut Denis Poisson, D\'epartement de Math\'ematiques, Universit\'e d'Orl\'eans, 45067 Orl\'eans Cedex~2, France}
\email{aline.bonami@gmail.com}

\maketitle

\begin{abstract}
This contribution starts with an exchange between us on the way we met Guido and he influenced our mathematical lives. Then it  is mainly a survey paper that illustrates this influence by describing different topics and their subsequent evolution after his seminal papers and courses.  Our main thread is the notion of a space of homogeneous type. In the second section we describe how it became central in pluricomplex analysis and consider particularly the existence  of weak factorization for spaces of holomorphic functions. In the last section, one revisits the construction of a basis of wavelets in a space of homogeneous type and the way it allows a Littlewood-Paley analysis.
\end{abstract}

\section{Introduction}
The authors of this chapter belong to two different generations of harmonic analysts. {We both met Guido at a very early stage in our careers.} His influence and friendship  were important for both of { us}.  {We found it interesting to cross our memories of Guido.} In a first section, we  will tell how we met Guido and try to describe what has been his influence on our mathematical choices.
 This  will be done through a dialog between us. 

The rest of the paper will be more classical and look like a survey paper on particular points for which the influence of Guido was like a starting point for us.  We will  not be exhaustive on the choice of topics and will certainly not cite all the appropriate literature in the fields that we have chosen to describe in more detail. We apologize in advance for these choices, which are linked with our personal views and memories and may diverge from other choices and memories. 

The main thread of this paper is the notion of space of homogeneous type that was introduced by Raphy Coifman and Guido in the early seventies. In the second section, mainly written by Aline, we will speak of its influence in pluricomplex analysis. In the third part, mainly written by Pascal, we will see how it may lead to develop all necessary tools for Littlewood-Paley theory in a very general context.

\section{The past, the past...}

\noindent{\bf Pascal.} You start, of course!\\

\noindent{\bf Aline.} For me, it begins in the academic year 1970-1971, when Guido Weiss and Raphy  Coifman gave a course in Orsay, now University Paris-Saclay. Their collaboration had started a few years earlier, and at the time they were interested in developing multiplier theorems related to actions of non-commutative groups. I was in the audience, which was not so numerous: even if harmonic analysis was still one of the main topics in the maths department, there was not such a feverish atmosphere as three or four years before. Moreover, Jean-Pierre Kahane and Yves Meyer were mainly abroad. It was perhaps a chance for us. Guido and Raphy asked for volunteers to help them in writing notes, which I did in collaboration with Jean-Louis Clerc and Bernard Mischler. So I met Guido and Raphy every week. I can still see them at the blackboard, in their very different styles, Guido writing calmly with his magnificent writing long formulas for spherical harmonics, while  Raphy tried to convince us of the simplicity of notions.  The notes that followed, {\sl Analyse Harmonique Non-Commutative sur Certains Espaces Homogènes} \cite{CW} went much further than the initial program and, in a sense, were the first steps into a new paradigm, the world of {\sl spaces of homogeneous type.}\\

\noindent{\bf Pascal.} Which influence had this course on you?\\

\noindent{\bf Aline.} The scientific life in Orsay was so intense at that time that it is difficult to sort out influences. It was my third course on real analysis given by foreigners. Eli Stein had given one in 1967-1968, just when I started research, the course that led to his book {\em Singular integrals...}. I had adored this course and had thought at the time that it was really what I wanted to do. But how? It seemed so difficult and the school of Chicago seemed so much in advance! At the same time I followed a course of Yves { Meyer}, who started to ask me questions on multipliers that led to my thesis on what is now called  hypercontractivity. Working with Yves was a huge luck. At the time of the venue in Orsay of Guido and Raphy, I  defended my thesis and even worked during a few months in probability theory. In the mean time, I had also attended a course of N. Riviere, who died prematurely. He talked on singular integrals related to parabolic equations. After the course of Guido and Raphy, one saw precisely how  the two courses  I had followed before dealt with two examples of their {\sl Espaces de nature homog\`ene}. One had a very clear description of the geometric objects involved. In fact their course opened the door to  a lot of problems in the theory of singular integrals. Afterwards, it was obvious for me that this was the kind of mathematics I preferred.

It would be unfair to limit the influence of Guido and Raphy at Orsay that year to mathematics. From the beginning, and particularly at that time when relations between mathematicians were still  somewhat formal in France, we were all stuck by their  generosity and kindness, which they expressed directly to everybody. \\ 

\noindent{\bf  Pascal.} And what influence had this course on the harmonic analysis team in Orsay at the time? \\

\noindent{\bf Aline.} It had a direct influence on Jean-Louis Clerc, Noël Lohoue and me, both for the content of their course itself and for the way to do mathematics, to interact between us. As I said before, relationships at work in France were then not so natural. There was still the idea one should publish alone. We started to work together, the three of us, in line with the course, trying to extend to other contexts known properties in Fourier Analysis. Inspired by their course, we read together the book of Helgason {\sl differential geometry and symmetric manifolds.} And with Jean-Louis I studied Ces\`aro means for expansions in spherical harmonics. I went two months to Washington University (WU) during the spring 1972. Of course Guido and Raphy had a large influence on the harmonic analysts in Orsay, but a little later. In fact, Yves Meyer became a frequent visitor at WU and started to work with Raphy a little later, in 1974. Afterwards the interactions between WU and French harmonic analysts  increased.
\\

Again, friendship mixes with mathematics and it is impossible not to speak of this also. I visited  regularly WU after 1972, even if much more rapidly.  We met also frequently in Europe. To come back to the far past, I  particularly remember having been at WU in 1981, at the same time as Jean-Lin Journé. We were both invited in Guido's and Barbara's summer house for a kind of country party. Conversations, canoeing on the small lake,  Guido removing splinters in my hand with the same meticulousness he used to show in his classes..., there was such a charm, as in a film.  Guido's jokes made a lot to create a friendly atmosphere  around him. Nevertheless Guido's life had its share of mystery and suffering,  starting with his childhood in Italy and the intellectual heritage of his parents.  He  was also a man of conviction, deeply attached to the values of human solidarity. He was a model for many of us. 
\\

{\it Nearly twenty years later,  Pascal spends two  and a half academic years at WU in the context of his military service that all French men had to accomplish at that time. He is not the first student of Yves Meyer to do this: Jean-Lin Journé was there in 1981-82. Its influence on Jean-Lin may be measured when looking at the  notes of the course that he gave there and which were published under the title {Calder\'on-Zygmund Operators,
Pseudo-Differential Operators and the
Cauchy Integral of Calder\'on}.} \\

\noindent{\bf Aline.} It is time for me to ask you a question. Did you know  Guido before going to WU ? How was your encounter?\\

\noindent{\bf Pascal.} No, I did not know him, nor, I must say, had read his books or articles. I was working on wavelet theory at the time, during the years 1986-87 and 1987-88. The mathematical theory of wavelets was just starting and I focused on that. For my military service, one possibility was to do it abroad in a teaching or research institution. Yves suggested that I go to WU. He contacted Guido and I recall reading his enthusiastic reply (a mailed letter at the time) dated December 31, 1987. I was just stuck by the date. Then started some exchanges with Guido (by the new email tool!) in the spring to set up things for the next fall (and for arrangements with French military office). I also read his book with Eli Stein \cite{SW} before leaving. The encounter was made easy in all means by Guido's welcoming help. With common interests in tennis and bird-watching, it became a strong friendship.\\

\noindent{\bf Aline.} How was the scientific life in WU at the time? Were there seminars, for instance? How did you interact with others?\\

\noindent{\bf Pascal.} First, the mathematics department was a highly renowned place for harmonic analysis. Raphy was long gone, but Al Baernstein, Bj\"orn Dalhberg, Steven Krantz, Richard Rochberg, Mitch Taibleson were very active members. This department also attracted many doctoral students from abroad (Spain, Italy, Argentina, China..) and post-docs, some became good friends (Estella and Rodolfo Torres, Carlo Morpurgo, Marco Peloso, Maria-Jes\'us Carro...) and became influential mathematicians back home after having benefited from the formidable and stimulating scientific atmosphere. I was the only French person (the French system did not promote much foreign post-docs in  those days). Guido wanted me to learn new maths so he proposed to Maria-Jes\'us, who had done her thesis on interpolation theory, and me to work on some problems on multipliers around the Riesz means, transference and other things. This led to two articles. The activity in the department was intense: regular analysis seminar, colloquiums, lots of master courses, and meetings in Guido's office. I still have the notes of the course of Guido on transference, based on his work with Raphy.  All details were concisely given on the board: ``too many details'' was I thinking, but I realize now that this allowed the topic to profoundly print in my mind. He was always encouraging people so as to make them available for more (integrating the f+ according to the words of A. Zygmund).  I recall he brought me to visit Raphy at Yale, and this visit had tremendous impact for me, as an example   of how curiosity drives mathematics.  After my thesis,  I was trying to understand (solve, why not?) the Kato conjecture in any dimension using wavelet methods as Ph.~Tchamitchian had just reproved nicely the strongly linked T(b) Theorem of David-Journé-Semmes via regular adapted wavelets \cite{T}.   When wandering in the Yale math department, I found a table with help yourself free printed (the pdf's  did not yet exist and some references could be hard to find) articles that people were giving away. I discovered the article of McIntosh on bounded holomorphic calculus from the Proceedings of the Center for Mathematical Analysis of the Australian National University and also the article of Jacques-Louis Lions pursuing the interpolation theory of Kato for maximal accretive operators, which inspired me for later research. \\

\noindent{\bf Aline.} You have a paper with Guido on Wilson bases. How was it?\\

\noindent{\bf Pascal.} During my stay at WU, Guido was curiously NOT interested to work on wavelets, despite my attempts to suggest problems, like the one which I solved later showing that a mild sufficient condition on the wavelet generating an orthonormal basis suffices to conclude it arises from a multiresolution analysis, introducing at the same time a special series whose values are the dimensions of certain vector spaces, that was next baptised by Guido as the ``dimension function'' and which became a central tool as beautifully explained in Guido's reference book with Eugenio Hernandez  \cite{HeW}. That trip to Yale was when a new interest came in Guido's mind. Raphy explained to us the concept of local Fourier bases and Guido wanted to understand the calculations behind. We started to work out the details, trying in Guido's style to make them as conceptually simple and elegant as possible. The observation that these bases and the Meyer orthonormal wavelet basis were built from the same tricks fascinated him. Victor Wickerhauser joined the discussion and the article was ready for publication \cite{AWW}. This is how he plunged into the mathematical theory of wavelets, trying to exhibit its finest structures and as always in the most simple and accessible terms. This was when I left WU at the end of 1990...\\

\noindent{\bf Pascal.} You also have papers with Guido in the nineties. Can you say something on them?\\

\noindent{\bf Aline.} As you said, Guido started to work on the theory of wavelets in the nineties. Just before, he was particularly interested by characterizations  coming from Littlewood-Paley theory. I like very much his notes with Frazier and Jawerth in the  CBMS Regional Conference Series (1991), \cite{FJW}. He had all reasons to be  fascinated by wavelets and had questions related to constructions. Yes, we wrote two papers, the first one \cite{BSW} on   band-limited wavelets. Our joint work was pursued by Gustavo Garrig\'os who  was then his PhD student. Gustavo came later to Orl\'eans as a post-doc, which was really great for me. But then Gustavo and I worked on Bergman spaces in tube domains, where of course we used a lot his book with Eli Stein and the properties of Hardy spaces in tube domains given there.  \\

\noindent{\bf Aline.} From your point of view, what is the mathematical legacy of Guido? 
\\

\noindent{\bf Pascal.} Of course, the discoveries he made all along his career, some that we are going to elaborate on in this article afterwards, which paved the way to create a strong school in harmonic analysis and linked topics: interpolation, singular integrals, Hardy spaces, spaces of homogeneous type... In all of these, there was always the will to leave no black holes in his way of writing, going into every details even the simplest ones, to make the reading accessible. In that respect, he was particularly proud of the Chauvenet prize he obtained in 1967 \cite{W}, the highest American award for mathematical expository writing. For example, he told me that for his book with Eli Stein, the first to propose a synthetic exposition of the extension of harmonic analysis from one to several Euclidean dimensions that occurred in the sixties, he thought at length about notation  avoiding as much as possible using coordinates,  to present the computations in the most conceptual manner (and with almost no typos).  This notation is still accurate as of today. The legacy is then clear. He posed the solid bricks on which many new results can be elaborated. Still, despite there are now many references on each of these topics, my first reaction when I look for a reference or a proof is to browse his masterpieces as I am sure to find what I need (and sometimes understand some new points that I left aside on earlier reading). When I write research articles, I always have in mind the words of Guido: try to stay simple and reader friendly. \\

\noindent{\bf Pascal.} Same question.\\

\noindent{\bf Aline.}  
 I will only speak of my own heritage. Of course I mention first his books, primarily his book with Eli Stein and the Notes of Orsay, which played different roles in my personal Pantheon. I had the impression to have the Notes in mind, just as they were engraved during the courses of Guido and Raphy, while I came back to the book with Eli Stein to read it in all details each time I gave a course myself. The book with Eugenio reaches the same degree of perfection. May be you find elsewhere a better source for intuition, for example in Yves's books, but if you want to understand all details, the book of Eugenio and Guido makes the job remarkably. These are three books that I keep on hand in my library. Then there are the works of Guido that inspired me for long. I would mention first the paper in Annals with Raphy and Richard Rochberg on commutators and weak factorization. I have the impression to have turned around this paper part of my life. We will come back to this later. I would like also to mention some parts of his work that surprised me when I heard him speak of them and that I regret not to have understood  more deeply, not to have come back to them later. You know, this kind of sensation one has that there, there is something I do not understand enough, there is some new music with which I would like to familiarize myself. For example, I recall how Guido was excited with Block spaces, that is, spaces that are defined from an atomic decomposition, but with atoms that are not of mean $0,$ and for which one has almost everywhere convergence results, for instance. Of course, the subject has been deepened afterwards, in particular by Fernando Soria. 
 I have nevertheless the feeling I should have tried to understand better what was around. The same with interpolation. Being exhaustive is impossible. But I will try to show, in the next section, how the work of Guido influenced mine. And of course, as you said, there is a lot more in the legacy of Guido than a list of results. There is a way to do mathematics, to write and to talk, and to do much more than only mathematics, there is a way to be curious of life and others, there is a lot.


\bigskip

\section{Orsay's Course and excursion inside complex analysis}

\subsection{The Notes: a new paradigm, the one of space of homogeneous type.}
One of the main new notions of the course was the notion of {\em space of homogeneous type}. For  $X$ a topological space $X,$ one first calls  {\em  pseudo-distance (or quasi-distance depending on the authors}) a function $\rho:X\times X \mapsto [0,\infty)$  such that there exists $K\ge 1$ such that, for all $x,y, z$ in $X,$ 
\begin{enumerate}
    \item $\rho(x, y)=0 \qquad \Longleftrightarrow \qquad x=y$
    \item $\rho(x, y)=\rho(y,x)$
    \item  
    $\rho(x, z)\leq K(\rho(x,y)+\rho(y,z)).$
\end{enumerate}
The  pseudo-ball (or quasi-ball, or simply ball) $B(x,r)$ with center $x$ and radius $r$ is the set of points $y\in X$ such that $\rho(x, y)<r.$

We say that $X$ is a {\em space of homogeneous type} if $X$ is endowed by a pseudo-distance $\rho$ and by a Borel measure $\mu$ if 
\begin{enumerate}
    \item the balls $B$ centered at $x$ constitute a basis of neighborhoods of $x$
    \item the measure $\mu$ is doubling, that is, there exists a constant $A$ such that, for all $x\in X$ and $r>0,$
    $$\mu(B(x, 2r)\leq A\mu(B(x, r)).$$
\end{enumerate}
These definitions are sufficient to be able to define the maximal function of $f\in L^1(d\mu)$ by
$$Mf(x):=\sup_{r>0}\frac 1{\mu(B(x,r))}\int_{B(x,r)}|f(y)|d\mu(y),$$
to prove the maximal theorem in this context as a consequence of a covering lemma of Vitali type, to develop as well Whitney's covering lemmas and to make it possible to write a Calder\'on-Zygmund decomposition of a function $f\in L^1(d\mu).$ As in the theory of singular integrals in $\R^n,$ it is then possible to prove $L^p$ estimates for singular integrals once one assumes that one already knows them for a particular value $p_0>1.$ To go ahead and be able to have theorems $T(1)$ or $T(b)$ asks for new tools, and this is what we will tempt to describe later.

Calder\'on-Zygmund theory was developed at the same time by Coifman and de Guzman in \cite{CdG}.  This last one studied systematically covering lemmas in relation with differentiation \cite{dG}, which is  always a living field of study (see for instance \cite{rigot}).  

It should be emphasized that, in the seventies, one had in mind some generalizations of Calder\'on-Zygmund theory, but not all. For example, the idea that this could be used for discrete structures appeared much later. Nevertheless, there is already this idea in the Notes that the doubling property of the measure $\mu$    can be replaced by the following one, which is weaker and is expressed  only as a geometric property (named nowadays as the geometric doubling property): \\

\noindent $2'$. \hspace{0.2cm} There exists a constant $N$ such that, for all $x\in X$ and $r>0,$ the ball $B(x,2r)$ is covered by at most $N$ balls of radius $r.$ \\

This will be seen as the right point of view for some constructions later on in this text.
\medskip

There are also emerging notions in the course, which have been central in the further development of the theory. In particular, the definition of an atom of the Hardy space $H^1$ and the atomic decomposition of $H^1$ have only  been given by Coifman a little later \cite{Coifman}. But the notion of atoms can already be found in the Notes (page 89 for instance). It is considered there for Riemannian symmetric spaces, and only from examples that are built with  approximate identities. The property of singular integrals to transport them into molecules is already present. All this program is developed, this time with precise definitions  and statements, in the two fundamental texts \cite{CWbul, TaiW}.

\medskip

In the Notes, one of the first examples of space of homogeneous type which is  given is the unit sphere in $\C^n$ endowed with the Euclidean measure $d\sigma$ and the distance $d_b(z, \zeta):=|1-\langle z \, ,{\zeta}\rangle|.$ 
We note $\mathbb B_{n}$ the unit ball in $\C^n$, $\partial \mathbb B^n$ its boundary, that is, the unit sphere and $\langle\cdot, \cdot\rangle$ the Hermitian inner product on $\C^n$. The distance $d_b$ (the index $b$ stands for the boundary) is clearly invariant under the action of the group $SU(n)$ and the distance from the point $\mathbf 1=(1, 0, \cdots,0)$  satisfies the property
$$d(\mathbf 1, z)\simeq |\Im z_1| +|z'|^2, \qquad z=(z_1, z'),$$
meaning that the quotient of the two quantities lies between two uniform constants in a neighborhood of $\mathbf 1.$ Moreover, $\sigma (B(z, r))\simeq r^n$ and 
$$|z-\zeta|^2\lesssim d_b(z, \zeta)\lesssim |z-\zeta|.$$ This metric is   usually called as the Kor\'anyi's metric. It had in particular  already been used by A. Kor\'anyi  for generalizing Fatou Theorem of non tangential convergence for harmonic functions in the unit disc to harmonic functions in $\mathbb B_n$ with respect to the Laplace-Beltrami operator for the Bergman metric \cite{Kor}. Indeed, in the same way that the first question for spaces of homogeneous nature was about generalizations of the Lebesgue Differentiation theorem, the first question for holomorphic functions was about generalizations of Fatou Theorem.
\medskip

At this point, let us recall definitions before going further. For $\Omega$  a smooth bounded domain in  $\mathbb C^n$ and $0 < p < \infty$,  let $L^p(\Omega)$ denote the
Lebesgue space with respect to the Lebesgue measure $dV(z)$ and   $A^p(\Omega)$ be  the corresponding Bergman space, that is, the
closed subspace of $L^p(\Omega)$ consisting of  holomorphic functions.
For $0 < p < \infty$,  let $L^p(\partial\Omega)$ denote the
Lebesgue space  on $\partial\Omega$ with respect to the induced
surface measure $d\sigma$ and $H^p(\Omega)$  the Hardy space of holomorphic
functions on $\Omega$, with norm given by
$$\|f\|_{H^p}^p  := \sup_{0<\varepsilon<\varepsilon_0}
\int_{\Omega}|f(w-\varepsilon \nu_w)|^p d\sigma(w),$$
 where $\nu_w$ is the unit exterior normal vector at the point $w$, with $\nu_w$ seen as a vector in $\mathbb C^n$.
 The question of generalizing non tangential convergence to higher dimension was treated in the small book of E. Stein \cite {StHol} (1972), which has been a major source of inspiration for its adaptation of potential theory to this context and for  its list of open questions. For such smooth domains, functions in $H^p(\partial\Omega)$ have admissible limits a.e.: the admissible regions, roughly speaking, are defined in the same way as non tangential regions in the unit disc, except that the ordinary distance on the circle is replaced by the Kor\'anyi   by the pseudo-distance on the boundary $\partial \Omega,$ that is,
 $$d_b(z, \zeta):= |\langle \nu_z, z-\zeta \rangle|+ |z-\zeta|^2.$$ 
  It is not symmetric, but $d_b(z, \zeta)\simeq d_b( \zeta, z)$ and can be symmetrized. These admissible regions are optimal for strictly pseudo-convex domains (that is, domains that are locally diffeomorphic to strictly convex domains), but not for  weakly pseudo-convex domains in which there   may be a larger swelling of admissible regions in complex tangential directions and a need for a definition that takes this into account. See for instance \cite{NSWbound}. Because of admissible limits one may consider $H^p(\Omega)$ as a subspace of $L^p(\partial \Omega).$ 
 \medskip
 
 The next question concerns Bergman and Szeg\"o projections, which we denote by $P_B$ and $P_S$ and which are, respectively, the orthogonal projection from $L^2(\Omega)$ to $A^2(\Omega)$ and the orthogonal projection from $L^2(\partial\Omega)$ to $H^2(\Omega)$ (identified with a subspace of $L^2(\partial\Omega)$). They are respectively given by the Bergman kernel $B_{\Omega}(z, \zeta)$ and the Szeg\"o kernel $S_{\Omega}(z, \zeta).$ For the unit ball they are given by
 $$S(z, \zeta)=\frac{c_n}{(1-\langle z \, ,{\zeta}\rangle)^n}, \qquad B(z, \zeta)=\frac{c_n}{(1-\langle z \, ,{\zeta}\rangle)^{n+1}}.$$
 Kor\'anyi and V\'agi, independently of Coifman and Weiss, developed around the same time a theory of singular integrals on homogeneous spaces, which led to the $L^p$ inequalities for the Szeg\"o projection of the unit ball \cite{KV}.  This one is now seen as an immediate example of singular integral on $\partial \mathbb B^n.$   Indeed, it  is elementary to see that
$$|S(\mathbf 1, z) -S(\mathbf 1, \zeta)|\lesssim \frac{d_b(z,\zeta)^{1/2}}{|1-\mathbf 1\cdot \zeta|^{n+1}},$$
which, by invariance by the action of $SU(n)$, gives the required inequality for  a singular integral on the space of homogeneous type given by $\partial \mathbb B^n,$ the Kor\'anyi distance and the Euclidean measure.

In order to obtain $L^p$ estimates for the Bergman projection, one may also be consider it as a singular integral for the distance that is given by the ordinary distance in the radial direction and by $d_b$ on the  boundary $\partial \mathbb B_n.$ But $L^p(\mathbb B_n)$ estimates can be obtained in a much simpler way by using Schur's lemma. These kinds of estimates have been first developed by Forelli and Rudin  (see \cite{FR}).

This was the beginning of a long story, which we evoke now. From 1970, a considerable work has been done on Hardy and Bergman spaces as well as Szeg\"o and Bergman projections in smooth bounded pseudo-convex domains of $\mathbb C^n.$  The link between the adequate pseudo-distance at the boundary and the two kernels, which is obvious for the unit ball, asks for long studies in general.
  Estimates on the Bergman and Szeg\"o kernels have been developed in the context of strictly pseudo-convex domains, then pseudo-convex domains of finite type in $\C^2,$ then convex domains of finite type by Fefferman \cite{Fef}, Nagel, Rosay, Stein and Wainger \cite{NRSW}, McNeal and Stein \cite {McNS} and many others.

The work of Henkin \cite{Hen} on strictly pseudo-convex domains play a fundamental role for having explicit representation formulas. They lead to the construction of  support functions $H\in \mathcal C^{\infty}(\Omega\times U)$ with $U$ a neighborhood of the boundary, such that $H(\cdot, \zeta)$ is holomorphic on $\Omega$ for all $\zeta\in U$ and
$$C^{-1}d(z, \zeta)\leq  |H(z, \zeta)| \leq Cd(z, \zeta).$$
Here $d$ is defined in a neighborhood of the boundary by
\begin{equation}\label{support}
    d(z, \zeta)=\delta (z) + \delta(\zeta) +d_b(\pi(z), \pi(\zeta))
\end{equation}
where $\delta(z)$ is the distance of $z$ to $\partial \Omega$ and $\pi(z)$ is the closest point of $\partial \Omega$, which is well-defined in $U$ for $U$ small enough and $d_b$ is the adequate distance on the boundary.  When $\Omega$ is the unit ball, one can take simply $H(z, \zeta)=1-\langle z \, ,{\zeta}\rangle.$ It is well adapted to Szeg\"o or Bergman spaces, which have only singularities  at the boundary. For convex domains of finite type (with the adequate pseudo-distance on the boundary), such a construction has been given by Diederich and Forn\ae ss \cite{DF}. We will use these support functions later on. Unlike what happens in the unit ball, they are not found in connection with the Bergman kernel.

\subsection{Weak factorization, commutators and Hankel operators}
In this subsection we will concentrate on one of the other fundamental papers of Guido, namely the one on commutators \cite{commutators}. It was written in 1975 with Richard Rochberg and Raphy and, again, had a large influence on problems that have been considered later. Some of them have been solved recently, some are still open. It is among the top cited papers of Guido and has been generalized in many directly during the nearly passed  fifty years. One can only find here a personal point of view, which leaves many developments aside. 

The first question that is treated in this paper is the characterization of $BMO(\mathbb R^n)$ as the class of functions $b$ for which the commutator $[M_b,R_j]$ are bounded in $L^2(\R^n).$ Here $M_b$ stands for the multiplication by $b$ and $R_j$ is the $j-$th Riesz transform, that is, the convolution by the distribution $c_n  p.v.\frac{x_j}{|x|^{n+1}}.$ This has been fundamental and, almost twenty years later, revealed itself to be one of the ways to understand the div-curl lemma in relation with compensated compactness \cite{CLMS}.

We will consider in more detail the second part of the article, which deals with holomorphic functions and factorization. Let us recall that in dimension $1$ every function $f\in H^p(\mathbb B_1)$ may be written as the product $gh,$ with $g\in H^{p_1}(\mathbb B_1), \;h\in H^{p_2}(\mathbb B_1)$ and $\frac 1p= \frac 1{p_1}+ \frac 1{p_2}$. This is false in higher dimension (see \cite{Gow}). In this second part  one finds for the first time the idea to replace factorization by weak factorization. Namely, they prove that any function $f\in H^1(\mathbb B_n)$ can be written as
$$f=\sum_{j=1}^\infty g_j h_j, \qquad \sum_j\|g_j\|_{H^p}\|g_j\|_{H^{p'}}\simeq \|f\|_{H^1}.$$
Their proof uses atomic decomposition of Hardy spaces $H^1(\mathbb B_n).$ Namely, a holomorphic atom $A$ is the Szeg\"o projection of an atom $a$, that is, $A=P_S a$ with $a$ a bounded function  of mean zero, which is supported in a pseudo-ball $B\subset \mathbb B_n,$ and such that $\|a\|_\infty\leq |B|^{-1}.$ They prove that $H^1$ has an atomic decomposition, that is, $f\in H^1(\mathbb B_n)$ can be written in terms of atoms $a_j$ related to balls $B_j=B(w_j, r_j)$
$$f=c+\sum_{j=1}^\infty \lambda_j a_j, \qquad \sum_j |\lambda_j|\simeq \|f\|_{H^1}.$$ It is then  sufficient to factorize an atom. At this point their proof may be in some way simplified and can be seen as  an easy consequence of the existence of a support function that satisfies \eqref{support}, which makes it easier to find  generalizations to other domains. One uses an idea that emerged later on: it is possible to ask for the atoms to have more moments that vanish, so that $|P_S a_j|$, which is not supported in $B_j$,  decreases  rapidly outside  $B_j$. Then, one just takes  $g_j= H(z, \zeta_j)^{-l}, $ where $\zeta_j$ is at distance $r_j$ of $\partial \mathbb B_n$ and  $\pi(z_j)= w_j$. So for $l>n$ the norm of $g_j$  in $H^p$ is easily estimated, as well as the one of   $H(z, \zeta_j)^{l} Pa$ in $H^{p'}$ when the atoms has sufficiently vanishing moments. 

This kind of proof  extends easily to other domains, see \cite{BPS}. It allows one to generalize  weak factorization to $H^p$ for $p<1$ but does not work when there is no atomic decomposition of $H^p$, that is, for $p>1$, even in the unit ball $\mathbb B^n$. Up to our knowledge, the existence of weak factorization of $H^p$ is still an open problem in dimension $n>1$: 
\smallskip

\noindent{ \em For $p>1$ and $\frac 1p= \frac 1 {p_1}+ \frac 1 {p_2},$ is there weak factorization for $H^p(\mathbb B^n)$ with products of functions in $H^{p_1}(\mathbb B^n) $ and $H^{p_2}(\mathbb B^n)$ when $n\geq 2$?}
\smallskip{}

 Another way to prove weak factorization consists in proving boundedness of the Hankel operator, which is the anti-linear operator on $H^p$ defined by $H_b(f)=P_S(b\overline f).$ More precisely, the fact that $H^1$ may be weakly factorized with products of functions in $H^p$ and $H^{p'}$ is equivalent to the fact that $H_b$ is bounded on $H^p$ if and only if $b$ belongs to the dual of $H^1.$ The same considerations prove that the open problem given above is equivalent to the fact that the Hankel operator $H_b$ maps $H^{p_1}$ into $H^{p_2'}$ if and only if $b$ belongs to $H^{p'}$.
 
  In $\R^n$,  such characterizations  have been proved recently by Hyt\"onen \cite{Hyt-com} for  commutators $[M_b, R_j]$. It is straightforward that  these commutators map $L^{p}$ into $L^{q}$ for $b\in L^r$, with $\frac 1r= \frac 1p +\frac 1q$. The challenge was to prove the converse, which is quite intricate. Unfortunately, this result, even if proved for the Szeg\"o projection on $\partial \mathbb B_n$, is not sufficient to find weak factorization for $H^p(\mathbb B_n)$ for $p>1$.

 Weak factorization for the Hardy space $H^1$ has also been studied in the polydisc, where there is no atomic decomposition. It was announced in a celebrated paper of Ferguson and Lacey but at this moment there is a gap in their proof as proved in \cite{HTV}. It is quite curious that the weak characterization of $H^1$, if true, implies, again on the bidisc, a weak factorization of $H^p$ by products of functions in $H^{2p}$ when $p<3/2,$ see \cite{BPSW}.
\medskip

The paper of Coifman,  Rochberg and Weiss finishes with the weak factorization of Bergman spaces. This is by now also completely classical and has been generalized in many contexts. The classical method, now, (see for instance \cite{Zhu}) goes through the atomic decomposition that has been developed by Coifman and Rochberg  \cite{CoifRoch} a little later: for $p>0$ given, every function in $f\in A^p(\mathbb B_n)$ may be written as
$$f(z)=\sum_j\lambda_j \frac {\omega_j(1-|w_j|)^N}{(1-\langle z \, ,{w_j}\rangle)^{n+N+1}}, \qquad \sum |\lambda_j|^p\simeq \|f\|_{A^p}^p,$$
where the points $z_j$ do not depend of the function $f$ but constitute what is called an $\eta-$net: they are at distance at most $\eta$ for the Bergman distance, but balls of radius $\eta/2$ centered at those points are disjoint while balls of radius $2r$ are almost disjoint.  The quantities that appear here, $\frac {\omega_j(1-|w_j|)^N}{(1-\langle z \, ,{w_j}\rangle)^{n+N+1}},$ are the values at $(z, w_j)$  of the Bergman kernel  when the Lebesgue measure $dV$ is replaced by the weighted measure $(1-|z|)^NdV(z)$. The constant $\omega_j$ is a normalization factor so that its $A^p-$norm is $1$. 
One has not only the atomic decomposition, but the weak factorization of Bergman spaces for $p\leq 1$. Generalizations are valid in the same context as for Hardy spaces, as well as estimates for related Hankel operators, now defined as $H_b(f)=P_B(b\overline f)$. But one can go much further in many directions.

First, the fact that such decompositions are valid for $p>1$ and a tricky use of Rademacher coefficients has allowed Pau and Zhao \cite{PZ} to get estimates with loss for the Hankel operator, and so to obtain weak factorization of Bergman spaces of $A^p(\mathbb B^n)$ for all values of $p$.

Secondly, the validity of atomic decomposition and weak factorization of Bergman spaces can be proved far beyond the contexts for which one can conclude for  Hardy spaces, as outlined in \cite{CoifRoch}: one can find such decompositions in tube domains over symmetric cones, for instance,  but with restrictions that are due to the lack of validity of $L^p$ inequalities for the Bergman projection for all $p>1$. In particular, if $\Omega=\R^n + i\Gamma,$ where $\Gamma$ is a symmetric cone (think of the  forward light cone $y_0>\sqrt { y_1^2+\cdots +y_n^2}$), results related to the Bergman space may be adapted (see \cite{NaSeh}), but with new constraints that are linked to the validity of $L^p$ inequalities for the Bergman projection. The constraints for these $L^p$ inequalities given in  \cite{BBGRS} are the best possible in the case of the forward light cone as a consequence of  the decoupling inequalities of Bourgain and Demeter \cite{BD}.

We are then far from the geometry of spaces of homogeneous type and singular integrals, but all this finds its origin in the seminal papers written from the years 1970 by Guido and those researchers that gathered around him.

\section{Littlewood-Paley decompositions and wavelets}

\subsection{The context}
  
The tools based on covering arguments that only need balls allow, in principle, to forget about algebraic contents to focus more on geometry or distributions of points. Still, familiar tool boxes such as convolution and  the Fourier transform need to be replaced by something else.

First one needs approximations to the identity in relation to scales replacing usual mollifiers. Second to develop a rich function space theory and the action of operators on them (such as Calder\'on-Zygmund operators), one needs Littlewood-Paley decompositions. It took some time to develop those tools in the optimal generality offered by the definition of a space of homogeneous type, that is, without extra hypotheses. The aim of this section is to explain the evolution of the ideas and a  solution. We use  material from   Auscher-Hyt\"onen \cite{AH} and also the introduction of  Han-Li-Ward \cite{HLW}. The latter reference contains a rich bibliography as well. 

Recall that the definition of a space of homogeneous type involves a set $X$  equipped with a quasi-distance $d$  and a positive doubling Borel measure $\mu$ (Borel measure is defined in several ways in the literature: either it is a measure on the Borel $\sigma$-algebra and one works only with Borel measurable functions or it is an outer measure for which Borel sets are measurable (in the sense of Carath\'eodory): in the latter case, we add  the regularity property  that every measurable set is contained in a Borel set with same measure). As said,  extra hypotheses had to be imposed in the proofs. However, it seemed that this would not limit the range of applicability thanks to an article of Macias-Segovia that describes the structure of spaces of homogeneous type \cite{MS}. This explains why removing these extra conditions remained unexplored for some time. 

A popular assumption is that the volume of a ball $B(x,r)=\{y\in X : d(x,y)<r\}$ is comparable to a power $r^\alpha$ of its radius, uniformly in its center (this is called the Ahlfors-David regularity property for the measure). Another one is that the volume of a ball $B(x,Cr)$ exceeds $(1+\eps)$ times that the volume of $B(x,r)$ for some constants $1<C<\infty$, $\eps>0$ (this is called the reverse doubling property).  One or the other are satisfied in many cases. Nevertheless,
if one wants one or the other  for all $0<r$ less than a fraction of the diameter of $X$, this excludes atomic measures, and for example   $\mathbb{Z}$ equipped with the induced absolute value and counting measure.  Considering that Hardy-Littlewood maximal theorem's original setup was on $p-$summable sequences, this would be paradoxical not to have a theory covering all situations.   

The balls $B(x,r)$ define the topology on $X$ but they might not be open sets nor even Borel sets when $d$ is a not a distance. This can be remedied by  replacing them by their interiors in the definition of doubling. Alternately, one result of Macias-Segovia \cite{MS} says one can change the quasi-distance to a metrically equivalent one (actually, a power of a genuine distance) for which the balls are open. So this is not a major issue as this only affects the   quasi-distance and the doubling constants. 

A more serious issue is again due to the difference between the metric case and the quasi-metric case:  a distance is always Lipschitz-continuous (with respect to itself) from the triangular inequality; it is not necessarily the case of a quasi-distance which might not even be H\"older-continuous. Another result of Macias-Segovia \cite{MS} says that one can change the balls so that it becomes true and even the measure becomes Ahlfors-David regular.  This is why most articles afterwards claimed that their results were proved in a "general" space of homogeneous type. But this is not quite correct. First, this change preserves the topology but it is not metrically equivalent. Next, it changes the classes of H\"older-continuous functions and of Calder\'on-Zygmund kernels. An example in \cite{HLW} shows that doing this change might not be desirable if the quasi-distance has a special interpretation in terms of an explicit operator as one may loose the connection. In complex analysis as above, the quasi-distance is defined according to the geometric context and one does not want to change it.   This is, however,  in this context (with a new quasi-distance from \cite{MS}) that the proof of the T(b) Theorem by David-Journé-Semmes was given \cite{DJS}. They were able to implement an unpublished method of Coifman to produce approximations to the identity and Littlewood-Paley decompositions. More precisely, the approximations to the identities are families of operators $(S_k)_{k\in \mathbb{Z}}$ ($k\ge k_0$ if $X$ is bounded) with kernels $S_k(x,y)$ enjoying the following properties for some $C<\infty$, $\eps>0$ and all $k$ and points 
\begin{align*}
    &S_k(x,y)=0 \ \mathrm{if}\ d'(x,y)\ge C2^{-k}, \  \mathrm{and}\ |S_k(x,y)|\le C2^k
    \\
    & |S_k(x,y)- S_k(x,y')| \le C2^{k(1+\eps)} d'(y,y')^\eps
    \\
    & |S_k(x,y)- S_k(x',y)| \le C2^{k(1+\eps)} d'(x,x')^\eps
    \\
    &
    \int_X S_k(x,y)\, d\mu(y)=1 = \int_X S_k (x,y)\, d\mu(x)   \end{align*}
Here the quasi-distance $d'$ is the one provided by the Macias-Segovia result. If $D_k:=S_{k+1}-S_k$, then for all $1<p<\infty$ and $f\in L^p(X,\mu)$ ($k\ge k_0$  and $\int f\, d\mu=0$ when $X$ is bounded) 
\begin{equation*}
    C^{-1}\|f\|_p \le \|(\sum_k |D_kf|^2)^{1/2}\|_p \le C\|f\|_p.\end{equation*}
    Hence the map $f\mapsto (D_kf)$ acts as a Littlewood-Paley decomposition. 
From sometime on, these methods have been streamlined and generalized (see \cite{HS, DH}) in the same context. 

The next development has been to work on spaces with the reverse doubling property.  This case is of much interest for example on Lie groups where the volumes of balls may obey two different power laws whether we assume small or large radii (and the exponents are often called local and global dimensions respectively). This covers also the case of some complete Riemannian manifolds. The reverse doubling condition is indeed a quantitative expression of connectedness of the space: annuli between two co-centered balls can never be empty and non-emptyness self improves to having mass at the correct scale. There, the technology of Coifman can be developed nicely (see Han-M\"uller-Yang, \cite{HMY1}) in a sense that estimates on the kernel family $S_k(x,y)$ go through with slight modification in their formulations.   A decisive inequality for the treatment of singular integrals in reverse doubling spaces is that 
\begin{equation}\label{fail}
    \sum_{k\ge 0} (\mu(B(x,C^k r))^{-1} \lesssim (\mu(B(x,r))^{-1}
\end{equation}
as the series is controlled in a geometric manner. 
Again, such conditions exclude spaces of homogeneous types of "discrete" nature and more. By the early 2010, about 40 years after the seminal definition of a space of homogeneous type, although the covering arguments and maximal theorem worked beautifully, it was still not known how to develop Littlewood-Paley theory in full generality  (again, despite what most authors wrote) and thus the function spaces and operator theory of singular integrals that go with.     We shall see that  the solution came from the need to work in absence of the doubling condition. Let us explain how.

In parallel to the development of "regular" decompositions, dyadic analysis was made possible thanks to the introduction of "dyadic cubes" on spaces of homogeneous type, first for Ahlfors-David regular measures by G. David \cite{D} and then by M. Christ \cite{Christ} in full generality. The idea is that  containment into a unique parent dyadic cube and  having a bounded number of children dyadic subcubes can be obtained from a carefully defined partial order on a family of selected points  with conditions only involving the quasi-distance. These points become the ``centers" of the cubes, which are  sets contained in balls and containing balls with those centers and radii comparable to their diameters. This clever construction allowed a local formulation of the T(b) theorem that was best designed toward an analytic capacity result in mind.

In $\mathbb{R}^n$, the dyadic cubes are linked to the Haar wavelets which form an orthonormal basis of $L^2( dx)$ where $dx$ is Lebesgue measure. The construction of Christ naturally leads to an orthonormal basis with respect to the doubling measure at hand, of course with non regular functions. The beautiful theory of orthonormal wavelets on $\mathbb{R}^n$, launched by Meyer \cite{Me} after the first unnoticed discoveries of Str\"omberg \cite{Str}, allows one to build orthonormal bases of localized and regular functions: not only they have decay away from a dyadic cube (or even have compact support) but also they enjoy H\"older regularity. Some can even be infinitely smooth. Explicitly or implicitly, Fourier transform is used. In fact, translation invariance and dilation invariance play a crucial role and indeed regular wavelet bases can be built on stratified nilpotent Lie groups  \cite{Lem2}. 

Coming back to Haar wavelets, it can be shown they provide unconditional bases of $L^p(dx)$  when $1<p<\infty$. The situation is different when $p=1$ and $p=\infty$. When $p=1$, they are unconditional families in the real Hardy space $H^1(\mathbb{R}^n)$ but span a strict subspace called the dyadic Hardy space $H^1_{dyad}(\mathbb{R}^n)$. It had become a challenge to prove that $H^1(\mathbb{R}^n)$ also admits an unconditional basis. B. Maurey \cite{Maurey} solved the question by showing the two spaces are isomorphic, but without providing an explicit basis.   Str\"omberg's orthonormal basis of spline functions (that turn to fit the wavelet theory and were rediscovered by Battle and Lemarié \cite{Bat, Lem1}) precisely yields an explicit construction. In the dual range, $BMO(\mathbb{R}^n)$, the dual of $H^1(\mathbb{R}^n)$ from the theorem of C. Fefferman \cite{Fef2}, was known to be strictly contained in its dyadic version $BMO_{dyad}(\mathbb{R}^n)$. It has to do with the fact that Euclidean distance $|x-y|$ and the distance  induced by the dyadic cubes (the diameter of the smallest dyadic cube containing $x$ and $y$) are not equivalent. In his celebrated book (\cite{Ga}, p.417), J. Garnett implicitly mentioned that classical $BMO(\mathbb{R})$ is the intersection of three of shifted copies of  $BMO_{dyad}(\mathbb{R})$. By shifted, it is meant that one considers  dyadic grids $2^k Q_a + \ell$, $k\in \mathbb{Z}, \ell \in \mathbb{Z}$, where $Q_a$ is the shifted interval $a+[0,1[$ by $a\in [0,1[$. More than an amusing observation, this shifting property can be thought as a random process among dyadic grids and this became a key tool to develop the singular integral theory on subsets $X$ of Euclidean space equipped with a Radon measure satisfying the power law upper estimate $\mu(B(x,r))\le C r^\beta$ for all $x\in X$ and $r>0$ as evidenced  in the work of Nazarov-Treil-Volberg, \cite{NTV}. 
The construction of  random families of dyadic grids on other contexts would then allow to develop singular integral theory there. But randomness cannot come from translation and dilation invariance anymore. A first construction of  T. Hyt\"onen and H. Martikainen on quasi-metric sets with the geometric doubling property \cite{HM} was designed to prove a T(b) theorem on such sets equipped with a measure satisfying the upper doubling condition, that is $\mu(B(x,r))\le \lambda(x,r)$ where $\lambda(x,2r)\le C\lambda(x,r)$ (a condition that covers both the doubling condition for the measure itself and the  power law upper estimate). Another possible  construction of random grids was proposed in \cite{NRV} and T. Hyt\"onen and A. Kairema \cite{HK}   streamlined the construction in \cite{HM}.  

\subsection{A regular wavelet basis after all}

As mentioned, while dyadic decompositions occurred in full generality, neither Littlewood-Paley decomposition, nor regular wavelet bases were available in 2011 in a "general" space of homogeneous spaces.  The existence of an orthonormal  basis  of regular wavelets was too much of a dream to have been conjectured anywhere. In \cite{DH}, which was one of the most advanced text on the topic at this time,  it was even written that they are "out of reach" and  frames (a substitute notion where redundancy is allowed)  were obtained but using again the Macias-Segovia trick. Still, when Tuomas Hyt\"onen came to Orsay in 2011, I (Pascal) asked him the problem. How to overcome the lack of any algebraic structure? How even can we make sure that for a quasi-distance that is not H\"older-continuous, one can obtain profusion (dense classes) of H\"older-continuous functions? Should we restrict to the metric case?    It turned out that the idea  of random  dyadic grids opened the door to a solution in \cite{AH} that required five steps. 
\begin{enumerate}
    \item The first step is the construction of appropriate random dyadic grids using {\bf nested} collections of separated and dense subsets of $X$ at each given scale.
    \item The second one is  the construction of {\bf spline} functions that are H\"older-continuous  of (some) order in $(0,1)$ depending on the space.
    \item In the third step, once the doubling measure is given, the spline functions provide a {\bf multiresolution analysis} on $L^2(\mu)$ (without the translation and dilation requirements, of course).
    \item In the fourth step,   {\bf construct the wavelets} by transposing Hilbertian algorithms of Y. Meyer \cite{M,Mey} for the construction of spline wavelets on $\mathbb{R}^n$ to spaces of homogeneous type.
    \item In the fifth step, the {\bf localisation} of the wavelets is obtained  by extending from the metric case  to the quasi-metric case existing lemmas due to Demko \cite{Dem} to estimate the entries of the inverse of a band matrix.
\end{enumerate}
 Steps 1 and 2 do not  require measure information but only the geometric doubling property. The spline functions have interesting properties of interpolation and reproduction which permit to construct the wavelets through a multiresolution structure.
 The H\"older regularity of the wavelets is the one of the spline functions. The construction of random grids was slightly modified later by Hyt\"onen-Tapiola \cite{HT} in the metric case to allow arbitrary H\"older exponent in $(0,1)$.

This in the end yields the following result (gathering both \cite{AH, HT})

\begin{theorem}\label{th:wavelets} Let $(X,d,\mu)$ be any space of homogeneous type and call $A_0$ its quasi-triangle constant. Let $a:=(1+2\log_2 A_0)^{-1}$ or $a:=1$ if $d$ is Lipschitz-continuous. There exist $\eta\in (0,1)$ depending only on $A_0$ (any $\eta\in (0,1)$ works when $A_0=1$),  constants $\delta>0$ and $\gamma>0$ depending only on $A_0$,  and  $0<C<\infty$ depending on $A_0$ and the doubling constant, one has the following. 

There are non-decreasing  sets $\mathscr{X}^k$, $k\in \Z$ (and $k\ge k_{0}$ if $X$ is bounded) of  $\delta^{k}$-separated and $2A_0\delta^k$-dense points in $X$, such that  if $\mathscr{Y}^k:=\mathscr{X}^{k+1}\setminus\mathscr{X}^k$, labelled  as $\mathscr{Y}^k=\{y^k_\beta\}_{\beta}$, 
there exists an orthonormal basis of real-valued functions $\psi^k_{\beta}$, $k\in \Z$ (and $k\ge k_{0}$ if $X$ is bounded), corresponding to each $y^k_{\beta}\in\mathscr{Y}^k$,  of $L^2(\mu)$  (or the orthogonal space to constants  if $X$ is bounded)  having exponential decay
 \begin{equation*}
  \abs{\psi^k_{\beta}(x)}\leq \frac{C}{\sqrt{ \mu(B(y^k_{\beta},\delta ^k))}} \exp\big( -\gamma (\delta^{-k}{d(y^{k}_{\beta},x)})^{a}\, \big)
,
\end{equation*}
H\"older-regularity
\begin{equation*}
  \abs{\psi^k_{\beta}(x)-\psi^k_{\beta}(y)}\leq \frac{C}{\sqrt{ \mu(B(y^k_{\beta},\delta ^k))}}
    \Big(\frac{d(x,y)}{\delta^k}\Big)^{\eta} \exp\big( -\gamma (\delta^{-k}{d(y^{k}_{\beta},x)})^{a}\, \big)
,\qquad
    d(x,y)\leq\delta^k,
\end{equation*}
 and vanishing mean
 \begin{equation*}
 \int_{X} \psi^k_{\beta}(x)\, d\mu(x)=0, \quad k\in \Z, k \ge k_{0}, y^k_{\beta}\in\mathscr{Y}^k.
\end{equation*}
  \end{theorem}
 A subset $S$  of $X$ is $\eta$-separated if two distinct elements of $S$ have distance at least $\eta$, it is $\eta$-dense if any point in $X$ has distance to $S$ not exceeding $\eta$.   This theorem and its argument provide us with the desired Littlewood-Paley decomposition.  

  \begin{corollary} Fix $1<p<\infty$. The orthogonal projections $Q_k$ on the closed linear span $W_k$ of the $\psi^k_{\beta}$ form a Littlewood-Paley decomposition. More precisely, there is a constant $0<C<\infty$ such that  for all $f\in L^p(\mu)$, ($k\ge k_0$ and $\int_X f\, d\mu=0$ if $X$ is bounded) 
\begin{equation*}
C^{-1}\|f\|_p \le \|(\sum_{k\in \mathbb Z} |Q_kf|^2)^{1/2}\|_p \le C\|f\|_p.    
\end{equation*}
Moreover, the family  $\psi^k_{\beta}$ (together with the indicator of $X$ if $X$ is bounded) forms an unconditional basis of $L^p(\mu)$ and
$P_k=\sum_{j< k} Q_j $ ($j\ge k_0$ and add the orthogonal projection onto constants if $X$ is bounded) is a family of approximation to the identity with  kernels $P_{k}(x,y)$,  symmetric in $x,y$, 
having size estimates
\begin{equation}\label{Pk}|P_{k}(x,y)|\le  \frac{C}{\sqrt{\mu(B(x,\delta^k))\mu(B(y,\delta ^k))}} \exp\big(-\gamma (\delta ^{-k}{d(x,y)})^s \, \big),
\end{equation}
regularity estimates
$$|P_{k}(x,y)-P_{k}(x,y')|\le C \left({\frac{d(y,y')}{\delta ^{k}}}\right)^\eta
\left ( \frac{\exp\big(-\gamma (\delta ^{-k}{d(x,y)})^s \, \big)}{\sqrt{\mu(B(x,\delta^k))\mu(B(y,\delta ^k))}}+ \frac{\exp\big(-\gamma (\delta ^{-k}{d(x,y')})^s \, \big)}{\sqrt{\mu(B(x,\delta^k))\mu(B(y',\delta ^k))}}\right)$$
for some $C,\gamma,s$ and all $x,y,y'\in X$ and $k\in \Z$ (with $k\ge k_{0}$ if $X$ is bounded). 
Moreover
$$
\int_{X}P_{k}(x,y)d\mu(x)=1.
$$
The kernel $Q_{k}(x,y)$ of $Q_{k}$ 
has similar estimates with {$s$ changed to $a$} and the additional multiplicative exponential factor in the size estimate
\begin{equation}\label{holes}
  \exp\big(-\gamma ((\delta^{-k}d(x,\mathscr{Y}^k))^a+(\delta^{-k}d(y,\mathscr{Y}^k))^a)  \, \big),
\end{equation}
and the cancellation condition
$$
\int_{X}Q_{k}(x,y)d\mu(x)=0.
$$

\end{corollary}

The exponential factor in \eqref{holes} is the key new object quantifying the geometry of $X$. We shall come back to this. 

The decay of wavelets is only exponential and an open question is whether one can obtain bounded support up to relaxing orthonormality to bi-orthogonality if necessary. 


We shall now indicate the idea of the construction of the splines: the guiding principle is that  their interpolation properties, their reproducing properties and their regularity are obtained from  a probabilistic argument allowed from  the axioms in the  construction  of  random dyadic  cubes.

\subsection{The key ideas of the construction of the spline functions}

The classical piece-wise linear splines on $\R$ are generated by translation and dilation of the function
\begin{equation*}
  s(x)= x1_{(0,1]}(x)+(2-x)1_{(1,2)}(x).
\end{equation*}
Another formula is $s(x)= 1_{[0,1)}*1_{[0,1)}(x)$ and writing out this convolution yields 
\begin{equation*}
 s(x)= \int_0^1 1_{[0,1)}(x-u)\ud u = \int_0^1 1_{[u,u+1)}(x)\, du.
\end{equation*}  
 Random dyadic intervals of sidelength $1$, in the sense of Nazarov, Treil and Volberg {\cite[Sec.~9.1]{NTV}}, are precisely defined by translating the standard intervals $[k,k+1)$ by a random number $u\in[0,1)$. Thus the unit cube with left end at the origin, $[0,1)$, is translated to $[u,u+1)$. 
So one can  think of $s(x)$ as  an average of the indicators of random dyadic intervals, or in probabilistic notation
\begin{equation*}
 s(x)= \prob_u \big(x\in [u,u+1)\big),
\end{equation*}
where $\prob_u$ is a fancy notation for the uniform probability measure on $[0,1[$ (i.e. the Lebesgue measure).
This is the basic idea which will guide us in constructing splines for a quasi-metric space $X$ with the  the geometric doubling property.
We next provide some details and key points of the constructions, especially  of the splines. Proofs can be consulted in \cite{AH}.  

The {\bf reference dyadic points} are chosen with the additional requirement of {\bf nestedness} compared to the earlier references \cite{Christ,HK,HM}. 
For every $k\in\Z$, we choose a set of \emph{reference dyadic points} $x^k_{\alpha}$ as follows: For $k=0$, let $\mathscr{X}^0:=\{x^0_{\alpha}\}_{\alpha}$ be a maximal collection of $1$-separated points. Inductively, for $k\in\Z_+$, let $\mathscr{X}^k:=\{x^{k}_{\alpha}\}_{\alpha}\supseteq\mathscr{X}^{k-1}$ and $\mathscr{X}^{-k}:=\{x^{-k}_{\alpha}\}_{\alpha}\subseteq\mathscr{X}^{-(k-1)}$ be maximal $\delta^{k}$- and $\delta^{-k}$-separated collections in $X$ and in $\mathscr{X}^{-(k-1)}$, respectively. It is easy to show that 
for all $k\in\Z$ and $x\in X$, the reference dyadic points satisfy
\begin{equation*}
  d(x^k_{\alpha},x^k_{\beta})\geq\delta^k\quad\big(\alpha\neq\beta\big),\qquad d(x,\mathscr{X}^k)=\min_{\alpha}d(x,x^k_{\alpha})<2A_0\delta^k.
\end{equation*}
The second property is the $2A_0\delta^k$-density of the set $\mathscr{X}^k$.
Note that $\mathscr{X}^k\subseteq\mathscr{X}^{k+1}$, so that every $x^k_{\alpha}$ is also a point of the form $x^{k+1}_{\beta}$, and thus of all the finer levels. Note also that  $\mathscr{X}^k$ reduces to one point for $k\le k_0$  for some $k_0$ if $X$ is bounded. Still we can continue with the full chain of sets for $k\in \mathbb Z$. We denote $\mathscr{Y}^k:=\mathscr{X}^{k+1}\setminus\mathscr{X}^k$, and relabel these points as $\mathscr{Y}^k=\{y^k_\beta\}_{\beta}$. These points  will be the parameter set of our wavelets.

The {\bf reference partial order} among the pairs $(k,\alpha)$ is the same as in M. Christ \cite{Christ}. For $\delta>0$ small enough (the smallness is fixed at the end to make every step work, and depends only on $A_0$),
 each $(k+1,\beta)$ satisfies $(k+1,\beta)\leq(k,\alpha)$ for exactly one $(k,\alpha)$, in such a way that
\begin{equation}\label{eq:leqProperties}
 d(x^{k+1}_{\beta},x^k_{\alpha})<\frac{1}{2A_0}\delta^k\quad\Longrightarrow\quad(k+1,\beta)\leq(k,\alpha)
 \quad\Longrightarrow\quad d(x^{k+1}_{\beta},x^k_{\alpha})<2A_0\delta^k.
\end{equation}
The pairs $(k+1,\beta)$ with $(k+1,\beta)\leq(k,\alpha)$ are called the children of $(k,\alpha)$. Geometric doubling implies that their number is uniformly bounded. 

Randomness is made possible by distinguishing the reference points using {\bf labels} as in \cite{HK}. 
Points $(k,\alpha)$ and $(k,\beta)$ are called \emph{neighbours}, if they have children $(k+1,\gamma)\leq(k,\alpha)$ and $(k+1,\eta)\leq(k,\beta)$ such that $d(x^{k+1}_{\gamma},x^{k+1}_{\eta})<(2A_0)^{-1}\delta^k$. In this case, by the  quasi-$A_0$-triangle  inequality,
\begin{equation*}
\begin{split}
 d(x^k_{\alpha},x^k_{\beta})
 &\le A_0 d(x^k_{\alpha},x^{k+1}_{\gamma})
 +A_0^2 d(x^{k+1}_{\gamma},x^{k+1}_{\eta})
 +A_0^2 d(x^{k+1}_{\eta},x^k_{\beta}) \\
 &<2A_0^2\delta^k+\tfrac12 A_0\delta^k+2A_0^3\delta^k<5A_0^3\delta^k.
\end{split}
\end{equation*}
The number of neighbours that any point can have is also uniformly bounded.

Each pair $(k,\alpha)$ is equipped with two \emph{labels}. The primary label, denoted by $\operatorname{label}_1(k,\alpha)\in\{0,1,\ldots,L\}$, where $L$ is the maximal number of neighbours, is chosen in such a way that any two neighbours have a different label. The secondary label, denoted by $\operatorname{label}_2(k,\alpha)\in\{1,\ldots,M\}$, where $M$ is the maximal number of children, is chosen in such a way that no two children of the same parent have the same label.

As described above, we now want to perform a perturbation of the original $x^k_\alpha$ and $\leq$ so as to obtain a {\bf parametrized family} of similar objects, on which probabilistic statements can later be made. The parameter space is
\begin{equation*}
  \Omega=\Big(\{0,1,\ldots,L\}\times\{1,\ldots,M\}\Big)^{\Z},
\end{equation*}
with a typical point denoted by $\omega=(\omega_k)_{k\in\Z}$, where $\omega_k=(\ell_k,m_k)\in\{0,1,\ldots,L\}\times\{1,\ldots,M\}$.

The \emph{random new dyadic points} $z^k_{\alpha}(\omega)$ are defined by
\begin{equation*}
  z^k_{\alpha}(w):=\begin{cases}  x^{k+1}_{\beta} & \text{if }\operatorname{label}_1(k,\alpha)=\ell_k,\text{ and } (k+1,\beta)\leq(k,\alpha),\text{ and }\operatorname{label}_2(k+1,\beta)=m_k, \\
    x^k_{\alpha} &  \text{if }\operatorname{label}_1(k,\alpha)\neq\ell_k,\text{ or } \not\exists(k+1,\beta)\leq(k,\alpha)
    \text{ such that }\operatorname{label}_2(k+1,\beta)=m_k. \\
    \end{cases}
\end{equation*}
Clearly, this point depends  only on $\omega_k$ so we write $z^k_{\alpha}(\omega)=z^k_{\alpha}(\omega_k)=z^k_{\alpha}$ and use the most convenient notational choice.

Equip $\Omega$  with the natural probability measure $\mathbb{P}_\omega$, which makes all coordinates $\omega_k=(\ell_k,m_k)$  independent of each other and uniformly distributed over the finite set $\{0,1,\ldots,L\}\times\{1,\ldots,M\}$.
For $(k+1,\beta)\leq(k,\alpha)$ fixed, then
\begin{equation*}
  \mathbb{P}_\omega( z^k_\alpha(\omega)=x^{k+1}_\beta )\geq\frac{1}{(L+1)M}.
\end{equation*}
In other words, every old point on the level $k+1$ has a positive (and bounded from below) probability of being a new point on the level $k$.
The new points behave qualitatively like the reference points, only with slightly weaker constants in the separation and density:
\begin{equation*}
  d(z^k_{\alpha},z^k_{\beta})\geq\frac{1}{2A_0}\delta^k,\qquad\min_{\alpha}d(x,z^k_{\alpha})<4A_0^2\delta^k.
\end{equation*}
The \emph{new partial order} $\leq_{\omega}$, $\omega=(\omega_k)_{k\in\Z}$, is set up as follows, declaring that
\begin{equation}\label{eq:newLeq}
   (k+1,\beta)\leq_{\omega}(k,\alpha)\quad\overset{\operatorname{def}}{\Longleftrightarrow}\quad
   \begin{cases} d(x^{k+1}_{\beta},z^k_{\alpha}(\omega))<\tfrac14 A_0^{-2}\delta^k,\qquad \text{or} \\
     (k+1,\beta)\leq(k,\alpha) \text{ and } \not\exists\gamma: d(x^{k+1}_{\beta},z^k_{\gamma}(\omega))<\tfrac14 A_0^{-2}\delta^k.
   \end{cases}
\end{equation}
In other words, to find the new parent of $(k+1,\beta)$ for the new partial order $\leq_{\omega}$, we first check whether the reference point $x^{k+1}_{\beta}$ is close (within distance $\tfrac14 A_0^{-2}\delta^k$) to some new dyadic point $z^k_{\alpha}(\omega)$. If yes, then the corresponding $(k,\alpha)$ will be the new parent of $(k+1,\beta)$. If no such close point exists, then we simply use the original partial order $\leq$ to decide the parent of $(k+1,\beta)$. With this definition, one can observe that 
for any given $k,\alpha,\beta$, the truth or falsity of the relation $(k+1,\beta)\leq_{\omega}(k,\alpha)$ depends only on the component $\omega_k$ of $\omega$.
Thus, it is independent of any event that occurs at level $k+1$ and higher. This explicit definition of the original partial order $\leq$ and of $\leq_\omega$ only depends of  the geometric configuration of the points, and a condition analogous to \eqref{eq:leqProperties} is a \emph{consequence} of the definition:
\begin{equation*}
 d(z^{k+1}_{\beta},z^k_{\alpha})<\tfrac15 A_0^{-3}\delta^k\quad\Longrightarrow\quad(k+1,\beta)\leq_{\omega}(k,\alpha)
 \quad\Longrightarrow\quad d(z^{k+1}_{\beta},z^k_{\alpha})<5A_0^3\delta^k.
\end{equation*}
Iteration yields that for all $\ell\geq k$,
\begin{equation*}
 d(z^{\ell}_{\beta},z^k_{\alpha})<\tfrac16 A_0^{-4}\delta^k\quad\Longrightarrow\quad(\ell,\beta)\leq_{\omega}(k,\alpha)
 \quad\Longrightarrow\quad d(z^{\ell}_{\beta},z^k_{\alpha})<6A_0^4\delta^k.
\end{equation*}

With the auxiliary objects at hand, the \textbf{random dyadic cubes} are easy to define.
As in \cite{HK},   introduce three families of these cubes---the preliminary, the closed, and the open:
\begin{equation*}
  \hat{Q}^k_{\alpha}(\omega):=\{z^{\ell}_{\beta}(\omega):(\ell,\beta)\leq_{\omega}(k,\alpha)\},
\end{equation*}
\begin{equation*}
  \bar{Q}^k_{\alpha}(\omega):=\overline{\hat{Q}^k_{\alpha}(\omega)},\qquad
  \tilde{Q}^k_{\alpha}(\omega):=\operatorname{interior}\bar{Q}^k_{\alpha}(\omega).
\end{equation*}
Note that $\hat{Q}^k_{\alpha}(\omega)$, and hence $\bar{Q}^k_{\alpha}(\omega)$ and $\tilde{Q}^k_{\alpha}(\omega)$, only depends on $\omega_{\ell}$ for $\ell\geq k$. Of course, the name "cube" solely refers to the analogous situation in Euclidean space  without any other geometric meaning for such sets (that same construction in Euclidean  space would not necessarily end up with a cube).

The following theorem summarizes the above properties of the random dyadic cubes for a fixed parameter $\omega$ (that is, they separately form  dyadic grids), and supplements the key statement about their probabilistic behaviour under the random choice of $\omega\in\Omega$.

\begin{theorem}\label{thm:cubes}
For any fixed $\omega\in\Omega:=(\{0,1,\ldots,L\}\times\{1,\ldots,M\})^{\Z}$, the cubes satisfy the following relations of a dyadic grid: the covering properties
\begin{equation*}
  X=\bigcup_{\alpha}\bar{Q}^k_{\alpha}(\omega),\qquad
  \bar{Q}^k_{\alpha}(\omega)=\bigcup_{\beta:(k+1,\beta)\leq_{\omega}(k,\alpha)}\bar{Q}^{k+1}_{\beta}(\omega);
\end{equation*}
the mutual disjointness property
\begin{equation*}
  \bar{Q}^k_{\alpha}(\omega)\cap\tilde{Q}^k_{\beta}(\omega)=\varnothing\qquad(\alpha\neq\beta);
\end{equation*}
and the comparability with balls:
\begin{equation*}
  B(z^k_\alpha(\omega),\tfrac16 A_0^{-5}\delta^k)\subseteq\tilde{Q}^k_\alpha(\omega)
  \subseteq\bar{Q}^k_{\alpha}(\omega)\subseteq B(z^k_{\alpha}(\omega),6A_0^4\delta^k).
\end{equation*}
Moreover, when $\Omega$ is equipped with the natural probability measure $\prob_{\omega}$, we have for some $\eta\in(0,1]$ the small boundary layer property:
\begin{equation}\label{eq:smallBdry}
  \prob_{\omega}\Big(x\in\bigcup_{\alpha}\partial_{\eps}Q^k_{\alpha}(\omega)\Big)\leq C\eps^{\eta}
  \qquad\Big(\partial_{\eps}Q^k_{\alpha}(\omega)
  :=\{y\in\bar{Q}^k_{\alpha}(\omega):d(y,{}^c\tilde{Q}^k_{\alpha}(\omega))<\eps\delta^k\}\Big);
\end{equation}
and in particular the negligible boundary property:
\begin{equation*}
  \prob_{\omega}\Big(x\in\bigcup_{k,\alpha}\partial Q^k_{\alpha}(\omega)\Big)=0
  \qquad\Big(\partial Q^k_{\alpha}(\omega):=\bar{Q}^k_{\alpha}(\omega)\setminus \tilde{Q}^k_{\alpha}(\omega)\Big).
\end{equation*}
Eventually, the original dyadic point $x^k_\alpha$ may also be viewed as a `centre' of the random dyadic cubes $\bar{Q}^k_\alpha(\omega)$ for all $\omega\in\Omega$:
\begin{equation*}
    B(x^k_{\alpha},\tfrac{1}{8}A_0^{-3}\delta^k)\subseteq\bar{Q}^k_{\alpha}(\omega)\subseteq\bar{B}(x^k_{\alpha},8A_0^5\delta^k).
\end{equation*}
\end{theorem}

 The small boundary property may not be true for balls in general. It was already observed by G. David (without randomness) and used by M. Christ for singular integrals that  dyadic cubes should have this property and that this can be built-in from the partial order. It will be the key to obtain the H\"older-regularity exponent of our splines and wavelets. 
The {\bf construction of splines} on $X$, and the proof of their basic properties, is indeed an easy consequence of the preparations above.
For every $(k,\alpha)$, we define the \textbf{spline function}  $s^k_{\alpha}$ by
\begin{equation}\label{eq:defSpline}
  s^k_{\alpha}(x):=\prob_{\omega}\Big(x\in\bar{Q}^k_{\alpha}(\omega)\Big).
\end{equation}

\begin{theorem}
The splines \eqref{eq:defSpline} satisfy the following properties: bounded support
\begin{equation}\label{eq:splineSupport}
  1_{B(x^k_{\alpha},\tfrac{1}{8}A_0^{-3}\delta^k)}(x)\leq s^k_{\alpha}(x)\leq 1_{B(x^k_{\alpha},8A_0^5\delta^k)}(x);
\end{equation}
the interpolation and reproducing properties
\begin{equation}\label{eq:splineSum}
  s^k_\alpha(x^k_\beta)=\delta_{\alpha\beta},\qquad
  \sum_{\alpha}s^k_{\alpha}(x)=1,\qquad
  s^k_{\alpha}(x)
  =\sum_{\beta}p^k_{\alpha\beta}\cdot s^{k+1}_{\beta}(x)
\end{equation}
where $p^k_{\alpha\beta}=\prob_{\omega}\big((k+1,\beta)\leq_{\omega}(k,\alpha)\big)$ so that  $\sum_\beta p^k_{\alpha\beta}=1$;
and H\"older-continuity
\begin{equation*}
  \abs{s^k_{\alpha}(x)-s^k_{\alpha}(y)}
    \leq C\Big(\frac{d(x,y)}{\delta^k}\Big)^{\eta}.
\end{equation*}
\end{theorem}
Let us prove 
the \textbf{H\"older-continuity of the splines} to show how the probabilistic smallness of the boundary regions, as expressed by \eqref{eq:smallBdry}, comes into play.  Indeed,
\begin{equation*}
\begin{split}
  \abs{s^k_{\alpha}(x)-s^k_{\alpha}(y)}
  &=\Babs{\prob_{\omega}\Big(x\in\bar{Q}^k_{\alpha}(\omega)\Big)-\prob_{\omega}\Big(y\in\bar{Q}^k_{\alpha}(\omega)\Big)} \\
  &=\Babs{\int_{\Omega}\big(1_{\omega:x\in \bar{Q}^k_{\alpha}(\omega)}-1_{\omega:y\in \bar{Q}^k_{\alpha}(\omega)}\big)\ud\prob_{\omega}} \\
  &\leq \int_{\Omega}\big(1_{\omega:x\in \bar{Q}^k_{\alpha}(\omega),y\notin\bar{Q}^k_{\alpha}(\omega) }
     +1_{\omega:y\in \bar{Q}^k_{\alpha}(\omega),x\notin\bar{Q}^k_{\alpha}(\omega) }\big)\ud\prob_{\omega} \\
  &= \prob_{\omega}\Big(x\in \bar{Q}^k_{\alpha}(\omega),y\notin\bar{Q}^k_{\alpha}(\omega)\Big)
     +\prob_{\omega}\Big(y\in \bar{Q}^k_{\alpha}(\omega),x\notin\bar{Q}^k_{\alpha}(\omega) \Big) \\
  &\leq \prob_{\omega}\Big(x\in \partial_{d(x,y)\delta^{-k}}\bar{Q}^k_{\alpha}(\omega)\Big)
     +\prob_{\omega}\Big(y\in \partial_{d(x,y)\delta^{-k}}\bar{Q}^k_{\alpha}(\omega)\Big) \\
    &\leq C\Big(\frac{d(x,y)}{\delta^k}\Big)^{\eta}.\qedhere
\end{split}
\end{equation*}

So far, there is no measure on $X$ attached to this construction. 

\begin{proposition}\label{prop:density} Let $\mu$ be a non-trivial positive Borel measure on  a quasi-metric space with the geometric doubling property $X$, finite on bounded Borel sets.
Let $1\leq p<\infty$. Then, the linear span of the spline functions and in particular the space of H\"older-$\eta$-continuous functions with bounded support where $\eta$ is the H\"older exponent of the splines are dense in $L^p(\mu)$.
\end{proposition}

Recall that there is a precaution to take in the definition of a Borel measure. Of course, density of Lipschitz-continuous functions with bounded support is well-known in the metric case. Density of H\"older-continuous functions to some order with bounded support has been proved independently in \cite{MMMM} at about the same time by a different argument. It is not clear it had been proved before.

\subsection{The Meyer algorithm}

We continue the development of the spline theory in the presence of a non-trivial Borel measure $\mu$ on $(X,d)$ so that $(X,d,\mu)$ is a space of homogeneous type in the most general sense with the exception that we assume that balls are Borel sets. If not, we may replace them by their interior and this  only affects constants.

The splines provide a multiresolution analysis of $L^2(d\mu)$. This consists of all properties of a classical multiresolution analysis of Meyer \cite[Definition 2.1]{M2} and Mallat \cite{Mal}, to the extent that this definition is meaningful in a quasi-metric space context: the classical postulates dealing with translations and dilations, specific to the Euclidean space and the Lebesgue measure, are now meaningless.
Let $V_k$ be the closed linear span of $\{s^k_\alpha\}_\alpha$ {in $L^2(d\mu)$}. Then $V_k\subseteq V_{k+1}$, and
\begin{equation*}
  \overline{\bigcup_{k\in\Z}V_k}=L^2(d\mu),\qquad
  \bigcap_{k\in\Z}V_k=\begin{cases} \{0\}, & \text{if $X$ is unbounded}, \\ V_{k_0}=\{\mathrm{constants}\}, & \text{if $X$ is bounded}, \end{cases}
\end{equation*}
where $k_0$ is some integer.
Moreover, the functions $s^k_\alpha/\sqrt{\mu^k_{\alpha}}$ form a Riesz basis of $V_k$: for all sequences of numbers $\lambda_{\alpha}$, we have the two-sided estimate
\begin{equation*}
  \BNorm{\sum_{\alpha}\lambda_{\alpha}s^k_{\alpha}}{L^2(d\mu)}\eqsim\Big(\sum_{\alpha}\abs{\lambda_{\alpha}}^2 \mu^k_{\alpha}\Big)^{1/2},
\end{equation*}
with $\mu^k_{\alpha}:= \mu(B(x^k_{\alpha},\delta^k)).$
With this at hand, the algorithm is as follows.
\begin{enumerate}
   \item Define  the closed linear span $Y_k$  of the splines ${s}^{k+1}_\beta$ that vanish at all the points in $\mathscr{X}^{k}$. Then, $V_{k+1}= V_k\oplus Y_k$, where the sum is topological and let $S_k: V_{k+1}\to Y_k$ be the bounded projection corresponding to this decomposition.   
\item If $Q_k$ is the orthogonal projection   of  $V_{k+1}$ onto $W_k$, the orthogonal complement of $V_k$ in $V_{k+1}$, then $S_k: W_k\to Y_k$ is an isomorphism, with inverse $Q_k: Y_k\to W_k$.
\item Define the $L^\infty$-normalized pre-wavelets 
 $$\tilde \psi^k_{\beta}:=Q_{k}s^{k+1}_{\beta}$$
 where  the $s^{k+1}_{\beta}$ are the splines in $Y_k$ of item 2. 
 \item Orthonormalize the pre-wavelets.
\end{enumerate}
This algorithm comes with explicit constructions:
\begin{enumerate}
\item  Construct the bi-orthogonal Riesz basis $\{\tilde{s}^k_\alpha\}_\alpha$ to $\{{s}^k_\alpha\}_\alpha$  in $V_k$, where bi-orthogonal means 
\begin{equation}\label{eq:biortho}
\langle s^k_{\alpha}, \tilde s^k_{\beta}\rangle_{L^2(\mu)}= \delta _{\alpha,\beta}.
\end{equation}
One has
\begin{equation*}
\sqrt{\mu^k_{\alpha}}\tilde s^k_{\alpha}(x) = \sum_{\beta\in \mathscr{X}^{k}} M_{k}^{-1} (\alpha,\beta)  \frac{s^k_{\beta}(x)}{\sqrt{\mu^k_{\beta}}},
\end{equation*} 
with $M_{k}$ being the infinite matrix indexed by $\mathscr{X}^{k}\times \mathscr{X}^{k}$ (write by abuse  $\alpha \in \mathscr{X}^{k}$ for $ x_\alpha^k \in \mathscr{X}^{k}$)   with entries
$$
M_{k}(\alpha,\beta)=\frac{\langle s^k_{\alpha}, s^k_{\beta}\rangle_{L^2(\mu)}}{\sqrt{\mu^k_{\alpha}\mu^k_{\beta}}}.
$$
\item Write \begin{equation}\label{eq:Qkf}
 Q_{k}f= f- \sum_{\alpha\in \mathscr{X}^{k}} \langle f, \tilde s^k_{\alpha}\rangle_{L^2(\mu)} s^k_{\alpha}
 =f- \sum_{\alpha\in \mathscr{X}^{k}} \langle f,  s^k_{\alpha}\rangle_{L^2(d\mu)} \tilde{s}^k_{\alpha},
\end{equation}
 because the sum is the orthogonal projection of $f$ onto $V_{k}$. Compute the pre-wavelets.
\item Observe that one can use $\mathscr{Y}^{k}$ (write by abuse  $\beta \in \mathscr{Y}^{k}$ for $ y_\beta^k \in \mathscr{Y}^{k}$) as a label set for the pre-wavelets. Then
\begin{equation*}
  \psi^k_\alpha(x):=\sum_{\beta\in\mathscr{Y}^k}\widetilde{M}_k^{-1/2}(\alpha,\beta)\frac{\tilde\psi^k_\beta(x)}{\sqrt{\mu^{k+1}_\beta}},
\end{equation*} where $\widetilde{M}_k^{-1/2}(\alpha,\beta)$  are the entries of the inverse of the square root of 
 the positive self-adjoint matrix
 $$
\widetilde{M}_k(\alpha,\beta):=   \frac{\langle  \tilde \psi^k_{\alpha},  \tilde \psi^k_{\beta}\rangle_{L^2(\mu)}}{\sqrt{\mu_{\alpha}^{k+1}{\mu_{\beta}^{k+1}}}}
 $$
 indexed by  $\mathscr{Y}^{k}\times \mathscr{Y}^{k}$. Here, $\mu_{\alpha}^{k+1}=\mu(B(y^{k}_\alpha, \delta^{k+1}))\sim \mu(B(y^{k}_\alpha, \delta^{k}))$ and recall that the points $y^{k}_\alpha$ are the $x^{k+1}_\alpha$ that do not belong to $\mathscr{X}^{k}$.
\end{enumerate}
 
\subsection{Matrices with exponential decay on a quasi-metric discrete set}

The algorithm requires to invert matrices in certain classes. The following result is well-known in the metric case. The extension to the quasi-metric case comes with an interesting argument having the spirit of the Macias-Segovia work (see also \cite{PS} for a very neat argument). We give the proof to illustrate this.

\begin{lemma}\label{lemma2} Let $\Xi$ be a $1$-separated set in   a quasi-metric space $(X,d)$ with quasi-triangle constant $A_{0}$ having the geometric doubling property with constant $N$. Consider a matrix  $M=(M(\alpha,\beta))$ indexed by $\Xi\times \Xi$ such that there exists $c>0$ for which
$$
C=\sup_{(\alpha,\beta)} \exp\big(c d(\alpha,\beta)\big) |M(\alpha,\beta)| <\infty.
$$
Then $M$ is bounded on $\ell^2(\Xi)$. 
If $M$ is invertible, then there exists $c'>0$ such that  
$$
\sup_{(\alpha,\beta)} \exp\big(c' (d(\alpha,\beta))^s\, \big) |M^{-1}(\alpha,\beta)| <\infty
$$
with  $s=(1+\log_{2} A_{0})^{-1}$ or $s=1$ if $d$ is a Lipschitz-continuous quasi-distance. If, in addition, $M$ is positive self-adjoint, then the same conclusion holds for $M^{-1/2}$.
\end{lemma}

If $d$ is a genuine distance or a Lipschitz-continuous quasi-distance, then this follows from Theorem 5 in \cite{Lem2} and the remark that follows it, which extends an earlier result in \cite{Dem} for band-limited matrices, and $s=1$. The exponent $s$ is best possible in the general case. 

\begin{proof} That $M$ is bounded on $\ell^2(\Xi)$ follows from a simple application the Schur lemma and we skip it. 

We continue with the following observation. For $n\ge 1$, let $\kappa_{n}$ be the best  constant in the inequality
$$
d(\alpha_{0},\alpha_{n}) \le \kappa_{n }(d(\alpha_{0},\alpha_{1})+ d(\alpha_{1},\alpha_{2})+ \ldots +d(\alpha_{n-1},\alpha_{n}))$$
for every chain $(\alpha_{0},\alpha_{1}, \dots, \alpha_{n})$ of $n+1$ elements (not necessarily distinct) of $\Xi$. It is clear that $(\kappa_{n})$ is non-decreasing and $\kappa_{1}=1$, $\kappa_{2}\le A_{0}$. 
Moreover, using $d(\alpha_{0}, \alpha_{m+n}) \le A_{0}(d(\alpha_{0}, \alpha_{m})+ d(\alpha_{m}, \alpha_{m+n}))$, it follows that $\kappa_{m+n}\le A_{0}{\cdot\max\{\kappa_{m},\kappa_{n}\}}$. {Thus $\kappa_{2n}\le A_{0}\kappa_{n}$ and therefore, $\kappa_{2^j}\le A_{0}^j$. 
We conclude that $\kappa_{n}\le A_{0}^{1+\log_{2} n}= A_{0} n ^{\log_{2}A_{0}}$. Note also that if $d$ is $L$-Lipschitz, then $d(\alpha_0,\alpha_n)\leq d(\alpha_0,\alpha_{n-1})+Ld(\alpha_{n-1},\alpha_n)$, which by iteration gives $\kappa_n\leq L$ for all $n$.}

Now assume that $M$ is {positive} self-adjoint and invertible. In this case, one can write $M=h(I-A)$ with $h=(\|M\| + \|M^{-1}\|)/2$ a positive real number and $A$ a matrix with norm $r=(\|M\|- \|M^{-1}\|)/(\|M\| + \|M^{-1}\|)<1$. Moreover, the coefficients of $A$ have the same decay as those of $M$. Without loss of generality, we normalize $h=1$.  Develop 
$
(I-A)^{-1}$ in the Neumann series $\sum A^n$ and estimate the coefficients $A^n(\alpha,\beta)$, $n\ge 1$, $\alpha\ne \beta$, in two ways. First 
$|A^n(\alpha,\beta)|\le r^n$. Second, we have {
\begin{equation*}
\begin{split}
  |A^n(\alpha,\beta)|
  &\leq \sum_{(\alpha_{1}, \ldots, \alpha_{n-1})\in \Xi^{n-1}} C^n \exp \big(-c (d(\alpha,\alpha_{1})+ d(\alpha_{1},\alpha_{2})+ \ldots +d(\alpha_{n-1},\beta))\big) \\ 
  &\leq C^n\exp\big(-\frac{c}{2\kappa_n}d(\alpha,\beta)\big)\sum_{(\alpha_{1}, \ldots, \alpha_{n-1})\in \Xi^{n-1}} 
   \exp \big(-\frac{c}{2} (d(\alpha_{1},\alpha_{2})+ \ldots +d(\alpha_{n-1},\beta))\big) \\
  &\leq \tilde{C}^n\exp\big(-\frac{c}{2\kappa_n}d(\alpha,\beta)\big),
\end{split}
\end{equation*}
where we applied  $n-1$ times the inequality with $\varepsilon=c/2$,
 } 
$$
\sup_{\alpha \in X}{\exp\big( \varepsilon d(\alpha, \Xi)/A_{0}\big)} \sum_{\beta \in \Xi} \exp\big(- \varepsilon d(\alpha,\beta)\big) \le c(\varepsilon,A_{0},N)<\infty. 
$$
As $\kappa_{n}$ is non decreasing, we have for any  integer $n_{0}$ using the second estimate for $0\le n\le n_{0}$ and the first for $n>n_{0}$,
\begin{align*}
 |M^{-1}(\alpha,\beta)| &\le (n_{0}+1) \tilde C^{n_{0}} \exp\big(-\frac{c}{2\kappa_{n_{0}}} d(\alpha,\beta)\big) + r^{n_{0}+1} (1-r)^{-1} 
\end{align*}
{and $(n_{0}+1) \tilde C^{n_{0}}\leq D^{n_0}$} for some large constant $D>0$.
Choosing $n_{0}$ as the first integer such that the first term dominates, we see that $d(\alpha,\beta)\eqsim n_0\cdot\kappa_{n_0}\lesssim n_0^{1/s}$, where $s=1/(1+\log_2 A_0)$ (or $s=1$ if $d$ is Lipschitz-continuous, recalling that $\kappa_{n_0}\leq L$ in this case). Hence, for some constant $c'>0$, 
\begin{equation*}
  \abs{M^{-1}(\alpha,\beta)}\lesssim r^{n_0}
  \lesssim\exp\big(-c' d(\alpha,\beta)^s\, \big).
\end{equation*}

For $M^{-1/2}$, we use the power series $ (I-A)^{-1/2}=\sum c_{n}A^{n}$.  As $0\le c_{n}\lesssim n^{1/2}$, the argument is the same. 
Finally, if $M$ is not {positive} self-adjoint, then we use $M^{-1}=M^*(MM^*)^{-1}$, the remark that on a 1-separated set $d\ge d^s$
and that the product of matrices having the exponential decay with factor $c$ as in the statement has exponential decay with factor $c'\in (0,c/A_0)$. 
\end{proof}

The above lemma applies a first time to the matrix $M_k(\alpha,\beta)$ on the 1-separated set $\mathscr{X}^{k}$ for the renormalized quasi-distance $d_k$ given by $d_{k} (x,y):={\delta ^{-k}}{d(x, y)}$ on $X \times X$ and then a second time to the matrix $\widetilde M_k(\alpha,\beta)$ on the 1-separated set $\mathscr{Y}^{k}$ for the renormalized quasi-distance $d_{k+1}^s$. The constants are therefore uniform over all $k$'s. This provides us with the parameters $s$ and $a$ that one can see in the Corollary. This insures  the exponential decay in the statements and the preservation of the H\"older exponent throughout the construction.

\subsection{The completion of Guido and Raphy's program}

We come back to the importance of the distance to the point set $\mathscr{Y}^{k}$ in the estimates  for the kernel $Q_k(x,y)$ of the projection $Q_k$. The absence of the reverse doubling property is linked to the presence of holes at certain scales in $X$. It could be very well be that $\mathscr{Y}^{k}$ is empty (in which case  the distance to it is $-\infty$ by convention and there is no corresponding wavelets) or that $\mathscr{Y}^{k}$ contains no points in some regions  of $X$ where there is an absence of nearby points at scale exactly $\delta^{k+1}$. Therefore the distance to $\mathscr{Y}^{k}$ is larger than $\delta^{k+1}$ in those regions. In other words, the set $\mathscr{Y}^{k}$ may not be dense at the scale $\delta^{k+1}$. This can be quantified in some organisation for the growth of the measure of balls: there exists $\eps>0$ depending only on $A_0$ and the doubling constant such that for every $x\in X$ and $r>0$, there exists a decreasing sequence, finite or infinite, of integers $\{k_j\}_{j=0}^J$ such that $r\leq\delta^{k_0}<\delta^{k_1}<\ldots$ and
\begin{equation*}
  \mu(B(x,\delta^k))\gtrsim(1+\varepsilon)^j \mu(B(x,r))\quad\text{and}\quad
  d(x,\mathscr{Y}^k)+\delta^k\gtrsim\delta^{k_{j+1}}\quad\text{if }k_j\geq k>k_{j+1},
\end{equation*}
where we interpret $k_{J+1}:=-\infty$ if $J<\infty$. As one cannot tell the growth of this sequence $k_j$, there is no reason that \eqref{fail} still converges.
But with this organisation, one can show the substitute  inequality to \eqref{fail}:
For all $x\in X$ and all positive numbers $r,\nu, a, \gamma$, we have
\begin{equation*}
  \sum_{k:\delta^k\geq r}\mu(B(x,\delta^k))^{-\nu}\exp\big(-\gamma(\delta ^{-k}{d(x,\mathscr{Y}^k)})^a\,  \big)\lesssim \mu(B(x,r))^{-\nu}
\end{equation*}
with implicit constant independent of $x$ and $r$.
Such estimates apply to check Calder\'on-Zygmund kernel estimates.  For example, one can deduce
\begin{equation*}
  \sum_{k,\alpha}\abs{\psi^k_\alpha(x)\psi^k_\alpha(y)}
  \leq\frac{C}{\mu(B(x,d(x,y)))}.
\end{equation*}
The operators defined by $\psi^k_\alpha \mapsto \eps^k_\alpha \psi^k_\alpha$ where $\eps^k_\alpha=\pm 1$, have kernels formally given by 
$K(x,y)=\sum_{k,\alpha}{\eps^k_\alpha \psi^k_\alpha(x)\psi^k_\alpha(y)}$. They are bounded on $L^2(d\mu)$ and, from the estimate above and similar ones for the regularity of $K(x,y)$ from the regularity of the wavelets, are  Calder\'on-Zygmund operators. One applies the extrapolation result  on spaces of homogeneous type in full generality to obtain their  $L^p(d\mu)$, $1<p<\infty$, boundedness. As a consequence, the family of wavelets is an unconditional basis in $L^p(d\mu)$.  This was done in \cite{AH}. Also the wavelet decomposition of BMO via a Carleson measure characterisation and a proof of the T(1) theorem  were presented there. 

It was clear that this would also open the way to develop the functional spaces and operator theory in this context: Hardy, Besov, Triebel-Lizorkin spaces with various definitions and descriptions. Here is a non-exhaustive list of references where these aspects are treated \cite{CLL, FY, HLW, FYL, HLYY, HHLYY, HWYY}.
As said, the construction of regular wavelets in full generality offers the first regular Littlewood-Paley decomposition at the same time. The Corollary above also points out the difference in the kernel estimates for  the approximation to the identity ($P_k(x,y)$) and the dyadic block ($Q_k(x,y)$). In the first one the scale is given by $\mathscr{X}^{k}$ and size estimates \eqref{Pk} come with somehow classical (or expected) decay; in the second, in addition to the scale and the (expected) decay,  there is an extra term \eqref{holes} measuring the decay in the distance to $\mathscr{Y}^{k}=\mathscr{X}^{k+1}\setminus \mathscr{X}^{k}$. Nor this extra decay, neither the role of the difference set $\mathscr{Y}^{k}$ could  be detected under reverse doubling.
   This tells us how should the operators of a Littlewood-Paley decomposition look like and what should  the   Calder\'on reproducing formula be in full generality (see in particular \cite{HLW, HLYY, HWYY, HHLYY}). Its convergence  in an appropriate distributional sense can be established. Next,   the functional spaces can be defined with any such Littlewood-Paley decomposition (the orthogonality being not required)  and they are independent of this choice as it  gives an equivalent norm with any other choice. For example, this  shows the intrinsic character of the spaces; they are independent of the choice of the wavelet basis (Each time we select  nested sets $\mathscr{X}^{k}$ one can construct a wavelet basis. Hence, there are many different choices). Using another point of view (as the one described by Aline above), the set of points $(\delta^{k+1},y^k_\alpha) \in \mathbb{R_+}\times X$ furnishes a representation by a  hyperbolic network structure   in which one should work. From there, one can imagine developing tent space theory and more. 

In conclusion,  the intuition of Guido and Raphy behind their  definition of a space of homogeneous type, continuing a long line of development since Fourier series (see for example the preface by Y. Meyer in \cite{DH}) has proved to be fully correct, in the sense that  one can build {\bf all} the tool boxes  that harmonic analysts rely on in Euclidean spaces only using geometry and distribution of points. This also makes a  bridge between continuous and discrete situations where {\it ad hoc} methods were developed independently. Indeed, the same tools can be  {\bf universally} applied to familiar (for analysts) situations such as Ahlfors-David or reverse doubling measures but also to basic examples  that were ruled out by these hypotheses like discrete groups  $\Z$ or $\Z/p\Z$ and the multidimensional analogs (see, e.g., \cite{S}), the typical discrete metric structures arising in theoretical computer science (trees, graphs, or strings from a finite alphabet), $\Q_{p}$ from arithmetic \cite{AES}, or discrete approximations  (like those constructed in \cite{ACI}).   The missing ingredients came from the need to develop appropriate technology to work in non-doubling spaces. We would not call it the fall of the doubling condition as did J. Verdera in 2002 \cite{Verdera}, but rather the completion of the doubling theory.  

And last but not least,  one can do this using orthonormal wavelets (which are special molecules, another invention of Guido and Raphy). 

The completion of the program launched by Guido and Raphy makes the ground ready for more. For example, another challenge is: can one construct regular wavelets or Littlewood-Paley decompositions that would apply universally in those non-doubling spaces?
 
\bibliographystyle{acm}

\begin{thebibliography}{}

\end{thebibliography}


\begin{thebibliography}{1O}



\bibitem{ACI}
{\sc Aimar, H., Carena, M., and Iaffei, B.}
\newblock Discrete approximation of spaces of homogeneous type.
\newblock {\em J. Geom. Anal. 19}, 1 (2009), 1--18.

\bibitem{AES}
{\sc Albeverio, S., Evdokimov, S., and Skopina, M.}
\newblock {$p$}-adic multiresolution analysis and wavelet frames.
\newblock {\em J. Fourier Anal. Appl. 16}, 5 (2010), 693--714.

\bibitem{AH}
{\sc Auscher, P., and Hytönen, T.}
\newblock Orthonormal bases of regular wavelets in spaces of homogeneous type. 
\newblock {\em Appl. Comput. Harmon. Anal. 34} (2013), no. 2, 266--296. 

\bibitem{AWW}
{\sc Auscher, P., Weiss, G., and Wickerhauser,  V.}
\newblock Local sine and cosine bases of Coifman and Meyer and the construction of smooth wavelets. 
 \newblock {\em Wavelet Anal. Appl., 2,} 237--256, Academic Press, Boston, MA, 1992. 

\bibitem{Bat}
{\sc Battle, G.}
\newblock A block spin construction of ondelettes. {I}. {L}emari\'e functions.
\newblock {\em Comm. Math. Phys. 110}, 4 (1987), 601--615.

\bibitem{BBGRS} 
{\sc Békollé, D.  Bonami, A.,  Garrig\'os, G., Ricci, F.,  and Sehba, B.} 
\newblock Analytic Besov spaces and Hardy-type inequalities in tube domains over symmetric cones. 
\newblock {\em J. Reine Angew. Math. 647}, (2010), 25--56.

\bibitem{BSW}
{\sc Bonami, A.,  Soria, and Weiss, G.}
\newblock Band-limited wavelets.
\newblock {\em J. Geom. Anal. 3},  (2015), (1993), no. 6, 543--578.


\bibitem{BoNa}
{\sc Bonami, A. and Nana C.}
\newblock Some questions related to the Bergman
projection in symmetric domains.
\newblock {\em Adv. Pure Appl. Math.},  (2015), 191--197.

\bibitem{BPS}
{\sc Bonami, A. Peloso, M., and Symesak F.}
\newblock Factorization of the Hardy spaces and Hankel operators on convex domains in $\C^n$.
\newblock {\em J. Geom. Anal. 11}, no. 3, 363--397.


\bibitem{BPSW}
{\sc Bonami, A. Pott, S., Sehba B., and Wick, B.}
\newblock Factorization of the Hardy space $H^s(\mathbb T^2).$
\newblock Unpublished manuscript.




\bibitem{BD}
{\sc Bourgain, J. and Demeter, C.}
\newblock Decouplings for curves and hypersurfaces with nonzero Gaussian curvature.
\newblock {\em J. Anal. Math. 133},  (2017), 279--311.



\bibitem{CLL}
{\sc Chen, C., Li, J., and  Liao, F.}
\newblock Some function spaces via orthonormal bases on spaces of homogeneous type. 
\newblock {\em Abstr. Appl. Anal.} 2014 Art. ID 265378, 13 pp. 

\bibitem{Christ}
{\sc Christ, M.}
\newblock A {$T(b)$} theorem with remarks on analytic capacity and the {C}auchy
  integral.
\newblock {\em Colloq. Math. 60/61}, 2 (1990), 601--628.

\bibitem{Coifman}
{\sc Coifman, R.~R.}
\newblock A real variable characterization of $H^p$.
\newblock {\em Studia Math. 51} (1974), 269--274.


\bibitem{CdG}
{\sc Coifman, R.~R., and de Guzm\'an, M.}
\newblock Singular integrals and multipliers on homogeneous spaces.
\newblock {\em Rev. Un. Mat. Argentina 25} (1970/71), 137--143.



\bibitem{CLMS}
{\sc Coifman, R.~R., Lions, P.-L., Meyer, Y., and Semmes S.}
\newblock Compensated compactness and Hardy spaces.
\newblock {\em J. Math. Pures Appl. (9) 72}, (1993), no. 3, 247--286.

\bibitem{CoifRoch}
{\sc Coifman, R.~R., and  Rochberg, R.}
\newblock Representation theorems for holomorphic and harmonic functions in $L^p.$
\newblock {\em Astérisque 77}, (1980) 11--66.

\bibitem{commutators}
{\sc Coifman, R.~R., Rochberg, R., and Weiss, G.}
\newblock {Factorization theorems for Hardy spaces in several variables.}.
  \newblock {\'E}tude de certaines int{\'e}grales singuli{\`e}res.
\newblock {\em Ann. of Math. (2) 103}
     (1976), no. 3, 611--635.
  
  \bibitem{CW}
{\sc Coifman, R.~R., and Weiss, G.}
\newblock {\em Analyse harmonique non-commutative sur certains espaces
  homog\`enes}.
  \newblock {\'E}tude de certaines int{\'e}grales singuli{\`e}res.
\newblock Lecture Notes in Mathematics, Vol. 242. Springer-Verlag, Berlin,
  1971.
 
  \bibitem{CWbul}
{\sc Coifman, R.~R., and Weiss, G.}
\newblock {Extensions of Hardy spaces and their use in analysis.}
\newblock {Bull. Amer. Math. Soc. 83 (1977), no. 4, 569--645.}




\bibitem{D}
{\sc David, G.}
\newblock {\em Wavelets and singular integrals on curves and surfaces},
  vol.~1465 of {\em Lecture Notes in Mathematics}.
\newblock Springer-Verlag, Berlin, 1991.

\bibitem{DJS}
{\sc David, G., Journ{\'e}, J.-L., and Semmes, S.}
\newblock Op\'erateurs de {C}alder\'on-{Z}ygmund, fonctions para-accr\'etives
  et interpolation.
\newblock {\em Rev. Mat. Iberoamericana 1}, 4 (1985), 1--56.

\bibitem{Dem}
{\sc Demko, S.}
\newblock Inverses of band matrices and local convergence of spline
  projections.
\newblock {\em SIAM J. Numer. Anal. 14}, 4 (1977), 616--619.

\bibitem{DH}
{\sc Deng, D., and Han, Y.}
\newblock {\em Harmonic analysis on spaces of homogeneous type}, vol.~1966 of
  {\em Lecture Notes in Mathematics}.
\newblock Springer-Verlag, Berlin, 2009.
\newblock With a preface by Yves Meyer.

\bibitem{DF}
{\sc Diederich, K., and  Forn\ae ss J. E.}
\newblock {Support functions for convex domains of finite type}.
\newblock {\em Math.
Z. 230},  (1999), 145--164.
  
 
\bibitem{Fef}
 {\sc Fefferman, C.}
 \newblock The Bergman kernel and biholomorphic mappings of pseudoconvex domains.
 \newblock {\em Invent. Math. 26} (1974),  1--65.

\bibitem{Fef2}
 {\sc Fefferman, C.}
\newblock Characterizations of bounded mean oscillation. \newblock
{\em Bull. Amer. Math. Soc. 77} (1971), 587--588.
 
 \bibitem{FR}
  {\sc Forelli, F.; Rudin, W.}
\newblock Projections on spaces of holomorphic functions in balls.
\newblock {\em Indiana Univ. Math. J. 24} (1974/75), 593--602.
  

   \bibitem{FY}
 {\sc Fu, X., and  Yang, D.}
 \newblock Wavelet characterizations of the atomic Hardy space H1 on spaces of homogeneous type. 
 \newblock {\em Appl. Comput. Harmon. Anal. 44} (2018), no. 1, 1--37.
 
 \bibitem{FL}
 {\sc  Ferguson, S. and Lacey, M.}
\newblock A characterization of product BMO by commutators
\newblock {\em Acta Math.  189} (2002) ,  143--160.

\bibitem{FJW}
  {\sc Frazier, M., Jawerth, B.,  and Weiss, G.}
  \newblock {\em Littlewood-Paley theory and the study of function spaces.}
  \newblock CBMS Regional Conference Series in Mathematics, 79.  American Mathematical Society, Providence, RI, 1991.



\bibitem{FYL}
 {\sc Fu, X., Yang, D., and Liang, Y.}
 \newblock Products of functions in $BMO(X)$ and $H^1_{at}(X)$ via wavelets over spaces of homogeneous type. 
\newblock {\em J. Fourier Anal. Appl. 23} (2017), no. 4, 919--990.

  \bibitem{Ga}
  {\sc Garnett, J.}
  \newblock {\em Bounded analytic functions.}
  \newblock Pure and Applied Mathematics, 96. Academic Press, Inc. [Harcourt Brace Jovanovich, Publishers], New York-London, 1981.
  
 \bibitem{Gow}
  {\sc Gowda, M. Seetharama}
  \newblock { Nonfactorization theorems in weighted Bergman and Hardy spaces on the unit ball of $\C^n$ ($n>1$).}
  \newblock  {\em Trans. Amer. Math. Soc. 277} (1983), no. 1, 203--212. 
  
 
  \bibitem{GP}
  {\sc Grellier, S., and Peloso, M. M.}
  \newblock { Nonfactorization theorems in weighted Bergman and Hardy spaces on the unit ball of $\C^n$ ($n>1$)..}
  \newblock  {\em Illinois J. Math. 46} (2002), no. 1, 207--232. 
  
  \bibitem{dG}
  {\sc de Guzm\'an, M.}
  \newblock {\em Differentiation of integrals in $\R^n$.}
  \newblock { Lecture Notes in Mathematics, Vol. 481.} 
  Springer-Verlag, Berlin-New York, 1975. 
 
 
  
  

  \bibitem{HLW}
  {\sc Han, Y.~S., Li, J., and Ward, L.}
 \newblock Hardy space theory on spaces of homogeneous type via orthonormal wavelet bases. 
 \newblock{\em Appl. Comput. Harmon. Anal. 45} (2018), no. 1, 120--169.





\bibitem{HMY1}
{\sc Han, Y., M{\"u}ller, D., and Yang, D.}
\newblock A theory of {B}esov and {T}riebel-{L}izorkin spaces on metric measure
  spaces modeled on {C}arnot-{C}arath\'eodory spaces.
\newblock {\em Abstr. Appl. Anal.\/} (2008), Art. ID 893409, 250.

\bibitem{HS}
{\sc Han, Y.~S., and Sawyer, E.~T.}
\newblock Littlewood-{P}aley theory on spaces of homogeneous type and the
  classical function spaces.
\newblock {\em Mem. Amer. Math. Soc. 110}, 530 (1994), vi+126.

\bibitem{HW}
{\sc Han, Y.~S., and Weiss, G.}
\newblock Function spaces on spaces of homogeneous type.
\newblock In {\em Essays on {F}ourier analysis in honor of {E}lias {M}. {S}tein
  ({P}rinceton, {NJ}, 1991)}, vol.~42 of {\em Princeton Math. Ser.} Princeton
  Univ. Press, Princeton, NJ, 1995, pp.~211--224.

\bibitem{HLYY}
 {\sc He, Z., Liu, L., Yang, D., and Yuan, W.}  
 \newblock New Calder\'on reproducing formulae with exponential decay on spaces of homogeneous type. 
 \newblock {\em Sci. China Math. 62} (2019), no. 2, 283--350.

\bibitem{HWYY}
{\sc He, Z., Wang, F., Yang, D., and Yuan, W.}
 \newblock Wavelet characterization of Besov and Triebel-Lizorkin spaces on spaces of homogeneous type and its applications. 
 \newblock {\em Appl. Comput. Harmon. Anal. 54} (2021), 176--226.

\bibitem{HHLYY}
 {\sc He, Z., Han, Y. Liu, L., Yang, D., and Yuan, W.}
 \newblock A complete real-variable theory of Hardy spaces on spaces of homogeneous type. 
 \newblock {\em J. Fourier Anal. Appl. 25} (2019), no. 5, 2197--2267.
 
 \bibitem{Hen}
{\sc Henkin, G., and Leiterer, J.}
\newblock {\em Theory of functions on complex manifolds.}
\newblock Monographs in Mathematics, 79. Birkhäuser Verlag, Basel, 1984. 

\bibitem{HeW}
{\sc Hern\'andez, E., and Weiss, G.}
\newblock {\em A first course on wavelets. }
\newblock Studies in Advanced Mathematics. CRC Press, Boca Raton, FL, 1996.

\bibitem{HTV}
{\sc Holmes I., Treil S., and Volberg A.}
\newblock Dyadic bi-parameter repeated commutator and dyadic product BMO.
\newblock arXiv:2101.00763.


\bibitem{Hyt-com}
{\sc Hyt\"onen, T.}
\newblock The $L^p$ to $L^q$
boundedness of commutators with applications to the Jacobian operator. 
J. Math. Pures Appl. (9) 156, 351-391 (2021).
\newblock {\em J. Math. Pures Appl. (9) 156},  (2021), 351--391. 

\bibitem{HK}
{\sc Hyt\"onen, T., and Kairema, A.}
\newblock Systems of dyadic cubes in a doubling metric space.
\newblock {\em Colloq. Math. 126}, 1 (2012), 1--33.   

  




\bibitem{HM}
{\sc Hyt\"onen, T., and Martikainen, H.}
\newblock Non-homogeneous $Tb$ theorem and random dyadic cubes on metric
  measure spaces.
\newblock {\em J. Geom. Anal. 22 },  (2012), no. 4, 1071--1107.

\bibitem{HT}
{\sc Hyt\"onen, T., and Tapiola, O.}
\newblock Almost Lipschitz-continuous wavelets in metric spaces via a new randomization of dyadic cubes. 
\newblock {\em J. Approx. Theory 185 } (2014), 12--30. 

\bibitem{Kor}
{\sc Kor\'anyi, A.}
\newblock Harmonic functions on Hermitian hyperbolic space.
\newblock {\em Trans. Amer. Math. Soc. 135} (1969), 507--516.

\bibitem{KV}
{\sc Kor\'anyi, A. , and Vági, S.}
\newblock  Singular integrals on homogeneous spaces and some problems of classical analysis.
\newblock {\em Ann. Scuola Norm. Sup. Pisa Cl. Sci. (3) 25} (1972), 575--648.


\bibitem{Lem1}
{\sc Lemari{\'e}, P.~G.}
\newblock Ondelettes \`a localisation exponentielle.
\newblock {\em J. Math. Pures Appl. (9) 67}, 3 (1988), 227--236.

\bibitem{Lem2}
{\sc Lemari\'e, P.~G.}
\newblock {Base d'ondelettes sur les groupes de Lie stratifi\'es. (Basis of
  ondelettes on stratified Lie groups).}
\newblock {\em Bull. Soc. Math. Fr. 117}, 2 (1989), 211--232.

\bibitem{MS}
{\sc Mac{\'{\i}}as, R.~A., and Segovia, C.}
\newblock Lipschitz functions on spaces of homogeneous type.
\newblock {\em Adv. in Math. 33}, 3 (1979), 257--270.


\bibitem{McNS}
{\sc McNeal, D., and Stein, E. M.}
\newblock Mapping properties of the Bergman projection on convex domains of finite type.
\newblock {\em Duke Math. J. 73}, (1994), 177--199. 

\bibitem{Mal}
{\sc Mallat, S.~G.}
\newblock Multiresolution approximations and wavelet orthonormal bases of
  {$L^2({\bf R})$}.
\newblock {\em Trans. Amer. Math. Soc. 315}, 1 (1989), 69--87.

\bibitem{Maurey}
{\sc Maurey, B.}
\newblock Isomorphismes entre espaces $H^1$. 
\newblock
{\em Acta Math. 145} (1980), no. 1-2, 79–120. 


\bibitem{Me}
{\sc Meyer, Y.}
\newblock Principe d'incertitude, bases hilbertiennes et algèbres d'opérateurs. (French) [The uncertainty principle, Hilbert bases and operator algebras] 
\newblock {\em Séminaire Bourbaki, Vol. 1985/86. Astérisque No. 145-146} (1987), 4, 209--223.

\bibitem{Mey}
{\sc Meyer, Y.}
\newblock Ondelettes et fonctions splines.
\newblock In {\em S\'eminaire sur les \'equations aux d\'eriv\'ees partielles
  1986--1987}. \'Ecole Polytech., Palaiseau, 1987, pp.~Exp.\ No.\ VI, 18.

\bibitem{M}
{\sc Meyer, Y.}
\newblock Wavelets and operators.
\newblock In {\em Analysis at {U}rbana, {V}ol.\ {I} ({U}rbana, {IL},
  1986--1987)}, vol.~137 of {\em London Math. Soc. Lecture Note Ser.} Cambridge
  Univ. Press, Cambridge, 1989, pp.~256--365.

\bibitem{M2}
{\sc Meyer, Y.}
\newblock {\em Ondelettes et op\'erateurs. {II}}.
\newblock Actualit\'es Math\'ematiques. [Current Mathematical Topics]. Hermann,
  Paris, 1990.
\newblock Op{\'e}rateurs de Calder{\'o}n-Zygmund. [Calder{\'o}n-Zygmund
  operators].

\bibitem{MMMM}
{\sc Mitrea, D., Mitrea, I., Mitrea, M., and Monniaux, S.}
\newblock 
  Groupoid metrization theory. With applications to analysis on quasi-metric spaces and functional analysis. 
  \newblock {\em Applied and Numerical Harmonic Analysis.} Birkhäuser/Springer, New York, 2013. 


Nagel, Alexander (1-WI); Rosay, Jean-Pierre (1-WI); Stein, Elias M. (1-PRIN); Wainger, Stephen (1-WI)
Estimates for the Bergman and Szegő kernels in certain weakly pseudoconvex domains.
Bull. Amer. Math. Soc. (N.S.) 18 (1988), no. 1, 55–59.

 \bibitem{NRSW}
{\sc Nagel, A., Rosay J. P., Stein, E. M., and Wainger, S.}
\newblock Estimates for the Bergman and Szegő kernels in certain weakly pseudoconvex domains.
\newblock {\em Bull. Amer. Math. Soc. (N.S.) 18, } (1988), 55--59.


  \bibitem{NSWbound}
{\sc Nagel, A., Stein, E. M., and Wainger, S.}
\newblock Boundary behavior of functions holomorphic in domains of finite type.
\newblock {\em Proc. Nat. Acad. Sci. U.S.A. 78}, (1981), no. 11, part 1, 6596--6599.

\bibitem{NaSeh}
{\sc Nana, C. and Sehba, B. F.}
\newblock Toeplitz and Hankel operators from Bergman to analytic Besov spaces of tube domains over symmetric cones.
\newblock  {\em St. Petersbg. Math. J. 30}, (2019), 723--750.  




\bibitem{NRV}
{\sc Nazarov, F., Reznikov, A., and Volberg, A.}
\newblock The proof of $A_2$ conjecture in a geometrically doubling metric
  space.
\newblock  {\em Indiana Univ. Math. J. 62} (2013), no. 5, 1503--1533.


\bibitem{NTV}
{\sc Nazarov, F., Treil, S., and Volberg, A.}
\newblock The {$Tb$}-theorem on non-homogeneous spaces.
\newblock {\em Acta Math. 190}, 2 (2003), 151--239.

\bibitem{PS}
{\sc Paluszy{\'n}ski, M., and Stempak, K.}
\newblock On quasi-metric and metric spaces.
\newblock {\em Proc. Amer. Math. Soc. 137}, 12 (2009), 4307--4312.

\bibitem{PZ}
{\sc Pau, J., and Zhao, R.}
\newblock Weak factorization and Hankel forms for weighted Bergman spaces on the unit ball. 
\newblock {\em Math. Ann. 363},  (2015), 363--383.


\bibitem{rigot}
{\sc Rigot, S.} 
\newblock Differentiation of measures in metric spaces. 
\newblock New trends on analysis and geometry in metric spaces, 93--116,
Lecture Notes in Math., 2296,  Springer, Cham, (2022) 


\bibitem{S}
{\sc Steidl, G.}
\newblock Spline wavelets over {${\bf R},{\bf Z},{\bf R}/N{\bf Z}$}, and {${\bf
  Z}/N{\bf Z}$}.
\newblock In {\em Wavelets: theory, algorithms, and applications ({T}aormina,
  1993)}, vol.~5 of {\em Wavelet Anal. Appl.} Academic Press, San Diego, CA,
  1994, pp.~155--177.
  
  \bibitem{StHol}
  {\sc Stein, E. M.} 
\newblock {\em Boundary behavior of holomorphic functions of several complex variables.}
\newblock Mathematical Notes, No. 11. Princeton University Press, Princeton, N.J., 1972.

 \bibitem{SW}
  {\sc Stein, E. M., and Weiss, G.} 
\newblock {\em Introduction to Fourier analysis on Euclidean spaces.}
\newblock Princeton Mathematical Series, No. 32. Princeton University Press, Princeton, N.J., 1971.



\bibitem{Str}
{\sc Str{\"o}mberg, J.-O.}
\newblock A modified {F}ranklin system and higher-order spline systems on
  {${\bf R}^{n}$} as unconditional bases for {H}ardy spaces.
\newblock In {\em Conference on harmonic analysis in honor of {A}ntoni
  {Z}ygmund, {V}ol. {I}, {II} ({C}hicago, {I}ll., 1981)}, Wadsworth Math. Ser.
  Wadsworth, Belmont, CA, 1983, pp.~475--494.


\bibitem{TaiW}
{\sc Taibleson, M and Weiss, G.}
\newblock {The molecular characterization of certain Hardy spaces. Representation theorems for Hardy spaces pp. 67–149,}
\newblock {\em{Astérisque}, 77, Soc. Math. France, Paris, 1980.}

\bibitem{TW}
{\sc Taibleson, M and Weiss, G.}
\newblock {Spaces generated by blocks. }
\newblock {{\em Probability theory and harmonic analysis} (Cleveland, Ohio, 1983), 209–226, Monogr. Textbooks Pure Appl. Math., 98, Dekker, New York, 1986.}



\bibitem{T}
{\sc Tchamitchian, Ph.}
\newblock Ondelettes et intégrale de Cauchy sur les courbes lipschitziennes. (French) [Wavelets and the Cauchy integral on Lipschitz curves] 
\newblock {\em Ann. of Math. (2) 129} (1989), no. 3, 641–649.  


\bibitem{Verdera}
{\sc Verdera, J.}
\newblock The fall of the doubling condition in {C}alder\'on-{Z}ygmund theory.
\newblock In {\em Proceedings of the 6th {I}nternational {C}onference on
  {H}armonic {A}nalysis and {P}artial {D}ifferential {E}quations ({E}l
  {E}scorial, 2000)\/} (2002), no.~Vol. Extra, pp.~275--292.

\bibitem{W} 
{\sc Weiss, G.}
\newblock  Studies in Real and Complex Analysis, 
\newblock  {\em MAA Stud. Math.}, vol. 3, Math. Assoc. America, Washington, DC, 1965, pp. 124--178.

\bibitem{Zhu}
{\sc Zhu, K.}
\newblock {\em Spaces of holomorphic functions in the unit ball.} (English) 
Graduate Texts in Mathematics 226. ,  Springer, New York  (2005).


\end{thebibliography}

\end{document}